\documentclass[12pt]{amsart}
\usepackage{enumerate}
\usepackage{amsmath, amssymb, amsthm, amsfonts}
\usepackage{tikz}
\usepackage{geometry}
\usepackage[pdftex, bookmarks=true]{hyperref}
\usepackage{bbm} 
\usepackage{url}
\usepackage{cmap} 
\usepackage{verbatim}  
\usepackage{bm}  

\usepackage{marginnote}

\newtheorem{theorem}{Theorem}[section]

\newtheorem{corollary}[theorem]{Corollary}
\newtheorem{lemma}[theorem]{Lemma}

\newtheorem{claim}[theorem]{Claim}
\newtheorem*{claim*}{Claim}

\theoremstyle{definition}
\newtheorem{remark}[theorem]{Remark}
\newtheorem*{remark*}{Remark}

\newtheorem{definition}[theorem]{Definition}

\newcommand\ds{\displaystyle}

\newcommand\Eb{\mathbb{E}}  
\newcommand\Pb{\mathbb{P}}  

\newcommand\Rb{\mathbb{R}}

\newcommand\Nb{\mathbb{N}}
\newcommand\Gb{\mathbb{G}}
\newcommand\Gbconf{\mathbb{G}^{\textrm{conf}}}

\newcommand\Oc{\mathcal{O}}  
\newcommand\Pcedge{\mathcal{P}_\mathrm{edge}}

\newcommand\Xc{\mathcal{X}}
\newcommand\Yc{\mathcal{Y}}

\newcommand\dhat{\hat{d}}
\newcommand\fhat{\hat{f}}

\newcommand\comp{\mathrm{c}}

\newcommand\al{\alpha}
\newcommand\alh{\hat{\alpha}}
\newcommand\be{\beta}
\newcommand\ga{\gamma}
\newcommand\Ga{\Gamma}
\newcommand\de{\delta}
\newcommand\De{\Delta}
\newcommand\eps{\varepsilon}
\newcommand\ka{\kappa}

\newcommand\alstar{\alpha^{\star}}
\newcommand\nustar{\nu^{\star}}

\newcommand\alfm{\alpha^{\mathrm{FM}}}
\newcommand\alsm{\alpha^{\mathrm{SM}}}
\newcommand\alaug{\alpha^{\mathrm{aug}}}

\newcommand\muedge{\mu_\mathrm{edge}}
\newcommand\muvtx{\mu_\mathrm{vtx}}

\newcommand\nuvtx{\nu_\mathrm{vtx}}

\newcommand\pivtx{\pi_\mathrm{vtx}}
\newcommand\althin{\al_{\mathrm{thin}}}
\newcommand\psimain{\psi^{\mathrm{main}}}
\newcommand\psirem{\psi^{\mathrm{rem}}}

\newcommand\sm{\setminus}

\DeclareMathOperator{\Aut}{Aut}

\newcommand{\defeq}{\mathrel{\vcenter{\baselineskip0.5ex \lineskiplimit0pt
                     \hbox{\scriptsize.}\hbox{\scriptsize.}}}%
                     =}

\newcommand\ind{\mathbbm{1}} 

\DeclareFontFamily{U}{matha}{\hyphenchar\font45}
\DeclareFontShape{U}{matha}{m}{n}{
  <-6> matha5 <6-7> matha6 <7-8> matha7
  <8-9> matha8 <9-10> matha9
  <10-12> matha10 <12-> matha12
  }{}

\DeclareSymbolFont{matha}{U}{matha}{m}{n}
\DeclareMathSymbol{\Lt}{3}{matha}{"CE}

\title[Boosted second moment method]
{Boosted second moment method\\in random regular graphs}

\author{Bal\'azs Gerencs\'er}
\address{HUN-REN Alfr\'ed R\'enyi Institute of Mathematics, Budapest, Hungary and ELTE E\"otv\"os Lor\'and University, Budapest, Hungary}
\email{gerencser.balazs@renyi.hu}
\thanks{The first author was supported by the National Research, Development and Innovation Office (NRDI; grant KKP 137490).}

\author{Viktor Harangi}
\address{HUN-REN Alfr\'ed R\'enyi Institute of Mathematics, Budapest, Hungary and ELTE E\"otv\"os Lor\'and University, Budapest, Hungary}
\email{harangi@renyi.hu}
\thanks{The second author was supported by the MTA-R\'enyi Counting in Sparse Graphs ``Momentum'' Research Group, by NRDI 
(grant KKP 138270), and by the Hungarian Academy of Sciences (J\'anos Bolyai Scholarship).}

\begin{document}

\begin{abstract}
Determining the asymptotic independence ratio of random regular graphs is a key challenge in the area of sparse random graphs. Due to the interpolation method, we have very good upper bounds at our disposal, which are actually known to be sharp for sufficiently large degrees. However, we are still in need of good explicit lower bounds for specific degrees.

The classical approach by Frieze and {\L}uczak achieves a lower bound by first applying the second moment method to sparse Erd\H{o}s--R\'enyi graphs, and then cleverly transitioning from that model to regular graphs. They obtain an asymptotic formula (as the degree tends to infinity) but no explicit lower bounds are derived for specific degrees.

In contrast, in this paper, we apply the second moment method directly to random regular graphs. This approach has a number of advantages. First, we can numerically compute good explicit lower bounds for any given degree $d$. Moreover, we can even boost this lower bound by arguing that the obtained independent set has a certain spatial Markov property. One can then exploit this property by making local modifications to the independent set, resulting in substantial improvements, and beating the previous best bounds for any $d \geq 10$. Finally, this method gives finer asymptotics as $d \to \infty$ than the original Frieze--{\L}uczak approach. 

Moreover, these results can be useful even beyond the scope of the independence ratio due to the fact that independent sets with the Markov property may be used to construct other objects in random regular graphs. To demonstrate this, we consider the problem of decomposing random regular graphs into stars.
\end{abstract}

\maketitle

\section{Introduction} \label{sec:intro}

This paper is concerned with random regular graphs: let $\Gb_{N,d}$ be a uniform random graph among all $d$-regular simple graphs on the vertex set $\{1,\ldots,N\}$. A standard Azuma-type argument shows that the independence ratio\footnote{The independence ratio is the size of the largest independent set divided by the total number of vertices.} of $\Gb_{N,d}$ is highly concentrated around its expectation for any fixed $d$ and $N$. A classic result \cite{bayati2013interpolation} shows that for each degree $d \geq 3$, there exists a constant $\alstar_d$ such that the independence ratio of the random $d$-regular graph $\Gb_{N,d}$ converges in probability to $\alstar_d$ as $N \to \infty$. Determining $\alstar_d$, or at least providing good bounds for it, is a major challenge of the area. 

In \cite{bollobas1981independence} Bollob\'as used random regular graphs to prove the existence of large-girth regular graphs with small independence number. His argument actually shows that $\alstar_d \leq \alfm_d$, where the first moment bound $\alfm_d$ is the unique root of a certain function $\varphi_d \colon (0,1/2) \to \Rb$; see \eqref{eq:phi}. Asymptotically, as $d \to \infty$, it is easy to show that $\alfm_d$ decays as 
\begin{equation} \label{eq:al_asymptotic}
\frac{2}{d} \bigg( 
\log d - \log \log d + 1 - \log 2 + o_d(1) \bigg) .
\end{equation}
We mention here that we will work out sharper estimates for $\alfm_d$ in the Appendix; see Theorem~\ref{thm:alfm_bnd}. 

In \cite{frieze1990independence,frieze1992independence} Frieze and {\L}uczak used a second moment argument to obtain a lower bound matching \eqref{eq:al_asymptotic}. In conclusion, $\alstar_d$ is asymptotically equal to \eqref{eq:al_asymptotic}. 

So we have a good understanding of the asymptotic behavior of $\alstar_d$, but is it possible to determine the precise value of $\alstar_d$ for a given $d$? Non-rigorous statistical physics methods \cite{barbier2013hardcore} suggested that, for $d \geq 20$, the precise value of $\alstar_d$ is given by a certain 1-RSB (1-step replica symmetry breaking) formula. In their seminal work \cite{ding2016maximum} Ding, Sun, and Sly confirmed this prediction for very large $d$.\footnote{There is a phase transition below degree $20$: for $d \leq 19$ one can prove upper bounds that are stronger than the 1-RSB formula \cite{harangi2023rsb}. These RSB bounds are likely to be extremely close to the true value.} Although there is no closed-form expression for this 1-RSB value, \cite[Theorem 3.1]{ding2016maximum} proves the following asymptotic comparison to the first moment bound:
\begin{equation} \label{eq:alstar_vs_alfm}
\alstar_d = \alfm_d - \left( \frac{2}{e} \, \frac{\log d}{d} \right)^2
\left( 1 + \Oc \left( \frac{\log \log d}{\log d} \right) \right) .
\end{equation}

Note, however, that neither the Frieze--{\L}uczak, nor the Ding--Sly--Sun result is helpful when one is interested in $\alstar_d$ for some specific degree, say, for $d=100$. What concrete lower bounds do we have at our disposal? Most results in this direction attained lower bounds by analyzing certain \emph{local processes or algorithms}.\footnote{Many results actually apply to the more general setting of graphs of (essentially) large girth.} A notable tool is Wormald's differential equation method.\footnote{This method was introduced in \cite{wormald1995differential}. See \cite{warnke2019wormalds} for a concise treatment.} The difficulty with this method is that it tends to be very challenging to carry out the required computer analysis for more sophisticated algorithms. Due to such difficulties, some approaches treat only the smallest cases $d=3$ or $d=3,4$. Currently, the best lower bounds are due to Csóka \cite{csoka2016independent} for $d=3,4$; to Duckworth--Zito \cite{duckworth2009large} for $d=5,6,7$; and to Wormald \cite{wormald1995differential} for $d \geq 8$, apart from $d=16$ for which degree a better bound was achieved via star decompositions in \cite{delcourt2025decomposing}. 

\subsection{Our bounds}
In this paper we revisit the second moment method, which was used in \cite{frieze1992independence} to attain the asymptotic formula \eqref{eq:al_asymptotic} for $\alstar_d$. In fact, \cite{frieze1992independence} built on an earlier work \cite{frieze1990independence}, which addressed the same problem for sparse Erdős--Rényi random graphs, and the second moment method was actually used in that earlier work. Our approach applies the second moment method directly to random regular graphs. The statement of the result involves a function $f_{d,\al} \colon [0,\al] \to \Rb$ with a rather complicated but explicit formula; see \eqref{eq:f_intro}. 
\begin{theorem} \label{thm:sm_intro}
Let $d\geq 3$ and $\al \in (0,1/2)$. Assume that the function $f_{d,\al}(\be)$, defined in \eqref{eq:f_intro}, has the following property:
\begin{equation} \label{eq:the_cond}
f_{d,\al} 
\text{ attains its maximum at } \al^2 .
\end{equation}
Then the independence ratio of $\Gb_{N,d}$ is at least $\al-o_N(1)$ with probability $1-o_N(1)$ as $N$ goes to $\infty$. In particular, $\alstar_d \geq \al$. 
\end{theorem}
For each degree we want to find the largest $\al$ for which this holds. 
\begin{definition} \label{def:alsm}
For an integer $d \geq 3$ let 
\[  \alsm_d \defeq 
\sup \big\{\al \, : \, \text{\eqref{eq:the_cond} holds} \big\}
= \sup \big\{\al \, : \, 
f_{d,\al}(\be) \leq f_{d,\al}(\al^2) 
\;\forall \be \in [0,\al] \big\}  .\]
We refer to $\alsm_d$ as the \emph{second moment (lower) bound}.
\end{definition}

By Theorem~\ref{thm:sm_intro} we have $\alstar_d \geq \alsm_d$. How sharp is this bound? We will give a very precise asymptotic estimate for $\alsm_d$ but first let us focus on the obtained bounds for specific degrees. 

We will see that verifying condition \eqref{eq:the_cond} essentially requires comparing two local maxima of $f_{d,\al}$. Therefore, it is straightforward to determine $\alsm_d$ with high precision using a computer for a specific degree $d$. It turns out that the lower bound $\alsm_d$ would already provide new results for $d \geq 20$, but we can actually boost it further with an additional step. 

This boosting is made possible by a strengthened version of Theorem~\ref{thm:sm_intro}, which guarantees the existence of an independent set (of density $\alsm_d$) satisfying a certain spatial Markov property (see Theorem~\ref{thm:Markovian_indep}). The argument relies on a refinement of the second moment method, using a result of Backhausz--Bordenave--Szegedy \cite{backhausz2022typicality} on \emph{typical processes} (see Lemma~\ref{lem:boosted_smm}). Exploiting the Markov property, we perform local improvements to the independent set and compute explicitly the resulting gain in density. In Section~\ref{sec:augment_indep} we present a simple procedure that yields the following augmented bound:
\begin{equation} \label{eq:al_aug}
\alaug_d \defeq \alsm_d + 
(1-\alsm_d) p^d\frac{1+\big( 1- p^{d-1} \big)^d}{2} 
\text{, where }
p = \frac{1-2\alsm_d}{1-\alsm_d} .
\end{equation}
The density increase ($\alaug_d$ vs $\alsm_d$) is substantial for degrees up to $d \approx 50$ and noticeable up to $d \approx 1000$. For larger degrees the gain is not significant: as we will see in Section~\ref{sec:asymptotic_gain}, asymptotically we have $\alaug_d - \alsm_d = d^{-2+o(1)}$, which does not get us considerably closer to (the conjectured value of) $\alstar_d$. 

Table~\ref{tab:bound_comparision} shows our bounds $\alsm_d$ and $\alaug_d$ compared to the best bounds previously known. We also included the RSB upper bounds in the table (which are conjectured to be sharp for $d \geq 20$). Our lower bounds are fairly close: for $d=500$ we have $0.0190$ (lower) versus $0.0193$ (upper), meaning that the ``optimality rate'' of our lower bound is $98.4\%$.
\begin{table}
\centering
\begin{tabular}{c||r|r|r|r}
 \textbf{degree} & best previous low.bnd  & 
 \begin{tabular}{@{}c@{}}\textbf{our lower bound} \\ $\alaug_d$ ($\alsm_d$) \end{tabular} 
 & best upp.bnd & opt.rate\\
 \hline
$10$  & \cite{wormald1995differential} \quad $0.2573$ 
& $0.258371$ \;($0.245521$) & $0.281105$ & $91.91\%$\\
$11$  & \cite{wormald1995differential} \quad $0.2447$ 
& $0.246775$ \;($0.234916$) & $0.268311$ & $91.97\%$\\
$12$  & \cite{wormald1995differential} \quad $0.2335$ 
& $0.236447$ \;($0.225479$) & $0.256856$ & $92.05\%$\\
$13$  & \cite{wormald1995differential} \quad $0.2234$ 
& $0.227167$ \;($0.216997$) & $0.246525$ & $92.14\%$\\
$14$  & \cite{wormald1995differential} \quad $0.2143$ 
& $0.218767$ \;($0.209312$) & $0.237148$ & $92.24\%$\\
$15$  & \cite{wormald1995differential} \quad $0.2061$ 
& $0.211113$ \;($0.202301$) & $0.228588$ & $92.35\%$\\
$16$  & \cite{delcourt2025decomposing} \quad $0.2000$ 
& $0.204101$ \;($0.195866$) & $0.220737$ & $92.46\%$\\
$17$  & \cite{wormald1995differential} \quad $0.1916$ 
& $0.197645$ \;($0.189931$) & $0.213502$ & $92.57\%$\\
$18$  & \cite{wormald1995differential} \quad $0.1852$ 
& $0.191674$ \;($0.184431$) & $0.206810$ & $92.68\%$\\
$19$  & \cite{wormald1995differential} \quad $0.1793$ 
& $0.186131$ \;($0.179314$) & $0.200597$ & $92.78\%$\\
\hline
$20$  & \cite{wormald1995differential} \quad $0.1738$ 
& $0.180967$ \;($0.174538$) & $0.194809$ & $92.89\%$\\
$50$  & \cite{wormald1995differential} \quad $0.0951$ 
& $0.104315$ \;($0.102334$) & $0.109797$ & $95.01\%$\\
$100$ & \cite{wormald1995differential} \quad $0.0572$ 
& $0.065030$ \;($0.064305$) & $0.067446$ & $96.42\%$\\
$200$ & \cite[Thm 4]{wormald1995differential} \quad $0.0260$ 
& $0.039010$ \;($0.038762$) & $0.040018$ & $97.48\%$\\
$500$ & \cite[Thm 4]{wormald1995differential} \quad $0.0123$ 
& $0.019005$ \;($0.018948$) & $0.019307$ & $98.44\%$\\
$1000$ & \cite[Thm 4]{wormald1995differential} \quad $0.0068$ 
& $0.010771$ \;($0.010753$) & $0.010890$ & $98.92\%$\\
$2000$ & \cite[Thm 4]{wormald1995differential} \quad $0.0037$ 
& $0.006014$ \;($0.006009$) & $0.006061$ & $99.25\%$\\
$5000$ & \cite[Thm 4]{wormald1995differential} \quad $0.0017$ 
& $0.002736$ \;($0.002735$) & $0.002750$ & $99.54\%$\\
$10000$ & \cite[Thm 4]{wormald1995differential} \quad $0.0009$ 
& $0.001493$ \;($0.001492$) & $0.001498$ & $99.68\%$
\end{tabular}
\caption{Bounds on the independence ratio $\alstar_d$. The last column (optimality rate) shows the relative size of our best lower bound ($\alaug_d$) compared to the upper bound.}
\label{tab:bound_comparision}
\end{table}

As we mentioned---even though our bounds do not come in the shape of a closed formula---it is relatively simple to obtain precise numerical values with a computer. Our SageMath code is able to quickly compute the bound for any given degree $d \leq 10^9$. We uploaded the code to the following GitHub repository.

\begin{center}
\url{https://github.com/harangi/indratio}
\end{center}

The reader may run our code on the following interactive website, which also contains a table of our bounds for degrees up to $200$.

\begin{center}
\url{https://www.renyi.hu/~harangi/ind_ratio.htm}
\end{center}

We point out that most papers in the literature only reported their bounds for a small number of degrees. In fact, some of them were able to carry out the required computer-assisted analysis only for the smallest degrees. Prior to this paper, in order to get an explicit lower bound on $\alstar_d$ for a degree $d>20$, one has had to use a weaker bound provided by an explicit formula given in \cite[Theorem 4]{wormald1995differential}. 

To complement these numerical results, we will prove the following concrete estimate for the second moment (lower) bound $\alsm_d$. In the spirit of \eqref{eq:alstar_vs_alfm}, we compare it to the first moment (upper) bound $\alfm_d$. (In turn, the reader may find accurate estimates for $\alfm_d$ in Section~\ref{sec:fm_bound}.)
\begin{theorem} \label{thm:vs_alfm}
Let 
\begin{equation} \label{eq:eps_d}
\eps_d \defeq \frac{4\sqrt{2}}{e} \frac{(\log d)^{1/2}}{d^{3/2}} .
\end{equation}
Then, for any degree $d \geq 3$, it holds that  
\begin{equation} \label{eq:alsm_vs_alfm}
\alfm_d - \eps_d 
\leq \alsm_d
\leq \alstar_d 
\leq \alfm_d .
\end{equation}
\end{theorem}
The proof will require very delicate estimates\footnote{The bound $\alsm_d \geq \alfm_d - \eps_d$ is very accurate: the inequality would fail to hold (for large $d$) if we replaced the multiplicative constant in $\eps_d$ with anything smaller.}, occupying a significant portion of the paper (Sections \ref{sec:estimates}--\ref{sec:technical}). 

Note that the Ding--Sly--Sun result \cite{ding2016maximum} determines $\alstar_d$ precisely \textbf{for very large} $d$, showing that $\alstar_d = \alfm_d- 1/d^{2-o(1)}$; see \eqref{eq:alstar_vs_alfm}. Our lower bound \eqref{eq:alsm_vs_alfm} has the advantage that it produces an explicit bound \textbf{for every} $d$, although with the somewhat weaker exponent  $3/2-o(1)$.

\subsection{Applications}

The significance of the existence of Markovian independent sets of density $\al=\alsm_d$ goes beyond obtaining the lower bound $\alstar_d \geq \alaug_d$. One can start from a Markovian independent set and try to create other objects, and hence proving the existence of more complex structures in random regular graphs. 

In this paper we present one such application, concerning so-called \emph{star decompositions}. In Section~\ref{sec:star_decomp} we demonstrate how our approach can be used to prove that the edge set of a random $d$-regular graph \emph{asymptotically almost surely (a.a.s)}\footnote{For a fixed $d$, we say that $\Gb_{N,d}$ a.a.s. has a property if it holds with probability $1-o_N(1)$.} can be decomposed into $k$-stars for certain pairs $d,k$. 

We also present a consequence of our numerical bounds on $\alstar_d$. There used to be a conjecture in sparse graph limit theory stating that so-called \emph{typical processes} are the same as the weak limits of \emph{factor of IID processes}. Loosely speaking, this says that ``most optimization problems over typical, sparse graphs can be solved by local algorithms'' \cite{rahman2017local}. The conjecture was refuted by a result of Gamarnik and Sudan \cite{gamarnik2017limits} showing that a typical independent set has larger density than a factor of IID independent set for sufficiently large $d \geq d_0$ (with no explicit $d_0$ derived). Our concrete bounds show that this statement already holds true starting at degree $d=403$. The details and precise statements can be found in the Appendix, see Section~\ref{sec:FIID_vs_typical}.

\subsection*{Organization of the paper}
\begin{itemize}
\item In Section~\ref{sec:counting} we provide a proof sketch for Theorem~\ref{thm:sm_intro}. Later we will rigorously prove a stronger version (Theorem~\ref{thm:Markovian_indep}) building on a result of Backhausz--Bordenave--Szegedy \cite{backhausz2022typicality} about typical processes. Although this section is not required for the rest of the paper, it may help the reader gain a clearer understanding, and it outlines a direct path to proving Theorem~\ref{thm:sm_intro}.
\item The next part of the paper contains our ``qualitative results''. Section~\ref{sec:typical} gathers the necessary facts about invariant and typical processes on the infinite tree. In Section~\ref{sec:sm_indep} we examine the second moment condition for independent sets, while in Section~\ref{sec:augment_indep} we discuss how the Markovian independent set can be augmented. 
\item Section~\ref{sec:star_decomp} is concerned with the application for star decompositions. 
\item The remainder of the paper (the ``quantitative part'') is dedicated to the proof of Theorem~\ref{thm:vs_alfm}. Section~\ref{sec:estimates} contains some general estimates. We analyze the function $f_{d,\al}$ and its local optima in Section~\ref{sec:local_maxima}, while Section~\ref{sec:second_max} focuses on the second local maximum and contains the main steps of the proof of Theorem~\ref{thm:vs_alfm}, with two technical arguments postponed until Section~\ref{sec:technical}. 
\item Finally, we included an Appendix containing sharp estimates for the first moment bound $\alfm_d$ (Section~\ref{sec:fm_bound}) and results about factor of IID independent sets (Section~\ref{sec:FIID_vs_typical}). 
\end{itemize}

\section{Counting independent sets in random regular graphs}
\label{sec:counting}

This section is somewhat independent of the rest of the paper, as its results are not needed in what follows, and stronger statements will be established later.

\subsection{Proof sketch of Theorem~\ref{thm:sm_intro}}
Our aim here is to outline a short, self-contained proof of our basic bound $\alstar_d \geq \alsm_d$ (Theorem~\ref{thm:sm_intro}). The argument relies only on standard results from the literature and makes direct use of the second moment method.

In a later section, we prove a stronger result (Theorem~\ref{thm:Markovian_indep}), where the second moment method is replaced by a theorem of Backhausz--Bordenave--Szegedy concerning typical processes with a spatial Markov property. That stronger result will imply the augmented bound $\alstar_d \geq \alaug_d$.

Recall that $\Gb_{N,d}$ denotes the uniform random graph among all $d$-regular simple graphs on the vertex set $\{1,\ldots,N\}$. This is closely related to the random multigraph $\Gbconf_{N,d}$, produced by the so-called \emph{configuration model}. Given $N$ vertices, each with $d$ ``half-edges'', the configuration model picks a random pairing of these $Nd$ half-edges, resulting in $Nd/2$ edges. The corresponding random graph $\Gbconf_{N,d}$, which is always $d$-regular, may have loops and multiple edges. A well-known fact is that if $\Gbconf_{N,d}$ is conditioned to be simple, then we get back $\Gb_{N,d}$. Moreover, for any fixed $d$, the probability that $\Gbconf_{N,d}$ is simple converges to a positive constant $p_d$ as $N \to \infty$. It follows that if $\Gbconf_{N,d}$ a.a.s.~has a certain property, then so does $\Gb_{N,d}$.

The advantage of working with $\Gbconf_{N,d}$ is that it makes counting fairly easy: for various graph structures or objects, one can write up their expected number in $\Gbconf_{N,d}$ using factorials (or multinomial coefficients). These quantities often have exponential growth (in $N$), the rate of which can be expressed using entropy.

In our setting, for fixed $d$ and $\al \in (0,1/2)$, let us define the random variable $Z_N$ as the number of independent sets of density $\al$ in $\Gbconf_{N,d}$, and consider its first and second moments: $\Eb Z_N$ and $\Eb Z_N^2$.

It is fairly easy to derive the exponential rate of $\Eb Z_N$: 
\begin{equation} \label{eq:phi}
\varphi_d(\al) \defeq 
h(\al) + \frac{d}{2} h(1-2\al) - (d-1) h(1-\al) , 
\end{equation}
where 
\[ h(x) = 
\begin{cases}
-x \log x & \text{if } x \in (0,1];\\
0 & \text{if } x=0 .
\end{cases} \]
More precisely, 
\[ \Eb Z_N = (Nd)^{\Oc(1)} 
\exp\big( N \varphi_d(\al) \big) . \]
The point is that the polynomial term $(Nd)^{\Oc(1)}$ is negligible in the sense that $\varphi_d(\al)<0$ implies that $\Eb Z_N \to 0$ as $N \to \infty$. Since the random variable is nonnegative and integer-valued, it follows that $\Pb( Z_N > 0 ) \leq \Eb Z_N$ also converges to $0$, that is, $\Gbconf_{N,d}$ (and hence $\Gb_{N,d}$) a.a.s. has no independent set of density $\al$. Recall that $\alfm_d$ is defined as the (unique) root of $\varphi_d(\al)$. Since $\varphi_d(\al)<0$ for every $\al>\alfm_d$, we can conclude that $\alstar_d \leq \alfm_d$. This is the first moment (upper) bound.

As for the second moment $\Eb Z_N^2$, we need to consider pairs of independent sets (both with density $\al$). We partition these pairs based on the intersection size. For $\be \in [0,\al]$ let $Z^{(2)}_{N,\be}$ be the number of pairs $(A_1,A_2)$ where both $A_1$ and $A_2$ are independent sets of $\Gbconf_{N,d}$ with density $\al$ such that their intersection $A_1 \cap A_2$ has density $\be$. It can be shown\footnote{The formulas for $\varphi_d$ and $f_{d,\al}$ stem from the entropy of certain discrete distributions. Section~\ref{sec:sm_indep} explains in detail which distributions need to be considered.} that
\[ \Eb Z^{(2)}_{N,\be} =
(Nd)^{\Oc(1)} 
\exp\big( N f_{d,\al}(\be) \big) , \]
where, for $\be \in [0,\al]$, let 
\begin{align}
\begin{split} \label{eq:f_intro}
\ga &\defeq -(1/2-\al)+\sqrt{(1/2-\al)^2+(\al-\be)^2} \text{, and} \\
f_{d,\al}(\be)  
&\defeq \frac{d}{2} \bigg( 2h(\be)+2h(\ga)
+4h(\al-\be-\ga) + h(1-4\al+2\be+2\ga) \bigg) \\
&\quad - (d-1) \bigg( h(\be)+ 2h(\al-\be) 
+ h(1-2\al+\be) \bigg) .
\end{split}
\end{align} 
Then the exponential rate of the number of pairs (with any intersection density $\be$) is 
\[ \varphi_d^{(2)}(\al) \defeq 
\max_{\be \in [0,\al]} f_{d,\al}(\be) .\]
By the Paley--Zygmund inequality\footnote{In this case this is simply the Cauchy--Schwarz inequality for the product $\ind_{\{Z_N>0\}} Z_N$.} 
we get 
\[ \Pb( Z_N > 0 ) 
\geq \frac{\big( \Eb Z_N \big)^2}{\Eb Z_N^2} 
= (Nd)^{\Oc(1)} \exp\big( N 
\big( 2\varphi_d(\al) - \varphi_d^{(2)}(\al) \big) \big) .\]
Note that we always have 
\[ \varphi_d^{(2)}(\al) \geq f_{d,\al}(\al^2) 
= 2 \varphi_d(\al) .\]
If $f_{d,\al}$ attains its maximum at $\al^2$, then we have equality, and obtain that 
\begin{equation} \label{eq:weak_ineq}
\Pb( Z_N > 0 ) \geq \frac{1}{(Nd)^K}  
\text{ for some constant } K. 
\end{equation} 
In other words, the independence ratio of $\Gbconf_{N,d}$ is at least $\al$ with probability $\geq (Nd)^{-K}$. Our goal, however, is to show that this holds a.a.s., that is, with probability $1-o_N(1)$, as $N \to \infty$. Next we show how $(Nd)^{-K}$ can be turned into $1-o_N(1)$ at the expense of slightly lowering $\al$.

As we mentioned in the introduction, the independence ratio of $\Gb_{N,d}$ is tightly concentrated around its expectation. We will use the following precise statement that readily follows from \cite[Theorem 2.19]{wormald1999survey}, which is a general Azuma-type concentration result. Let $i(G)$ denote the independence ratio of $G$; then 
\begin{equation} \label{eq:azuma}
\Pb\bigg( \big|i(\Gbconf_{N,d}) - \Eb \, i(\Gbconf_{N,d}) \big| \geq t \bigg) 
\leq 2 \exp\left( \frac{-Nt^2}{d} \right) .
\end{equation}
We set $t = t_N \defeq (\log N)/\sqrt{N}$ so that 
\[ \exp\left( \frac{-Nt^2}{d} \right) 
= \exp\left( \frac{-(\log N)^2}{d} \right) 
= o \left( \frac{1}{(Nd)^K} \right) 
\text{ for any fixed $d$ and $K>0$.} \]
Since the event $Z_N>0$ is equivalent to $i(\Gbconf_{N,d}) \geq \al$, from \eqref{eq:weak_ineq} and \eqref{eq:azuma} we obtain  
\[ \Eb \, i(\Gbconf_{N,d}) \geq \al - t_N .\]
Using \eqref{eq:azuma} again, we can conclude that with probability $1-o_N(1)$ we have 
\[ i(\Gbconf_{N,d}) 
\geq \Eb \, i(\Gbconf_{N,d}) - t_N 
\geq \al - 2t_N .\]
Since $t_N=o_N(1)$, the statement of Theorem~\ref{thm:sm_intro} follows.

\subsection{Labelings and entropy} \label{sec:labelings}

This section is primarily intended for readers who are unfamiliar with the language of graph limits used in the subsequent sections. Our aim is to explain the relationship between counting in $\Gbconf_{N,d}$ and entropic quantities in the context of graph labelings. 
It will shed some light on where the formulas for $\varphi_d$ and $f_{d,\al}$ come from.

Suppose that $\Xc$ is a finite set of labels. Given a $d$-regular graph $G$, an $\Xc$-labeling of $G$ is a mapping $\ell \colon V(G) \to \Xc$. If $G$ is finite, one can consider the local statistics of a labeling. For instance, we may define the \emph{edge statistics} of $\ell$ as the distribution $\pi$ of the pair $\big( \ell(u), \ell(v) \big)$, where $uv$ is a uniform random (directed) edge of $G$. Note that 
$\pi$ is an exchangeable distribution on $\Xc \times \Xc$ (i.e., it is invariant under exchanging the two coordinates). In particular, the two marginals of $\pi$ coincide; let us denote this marginal distribution by $\pivtx$. (Since $G$ is regular, $\pivtx$ is simply the distribution of $\ell(v)$ for a uniform random vertex $v$.) We denote the probabilities defining $\pi$ and $\pivtx$ as follows:
\begin{align*}
p_{x,y} &=p_{y,x}=\pi\big( \{(x,y)\} \big) ;\\
p_x &= \sum_{y\in \Xc} p_{x,y} 
= \pivtx\big( \{x\} \big) .
\end{align*}
The key quantities here are the \emph{Shannon entropy} of $\pi$ and $\pivtx$ defined as follows:
\begin{equation} \label{eq:Sh_entr}
H(\pi) = \sum_{x,y \in \Xc} h\big( p_{x,y} \big) 
\quad \text{and} \quad
H(\pivtx) = \sum_{x \in \Xc} h\big( p_x \big) .
\end{equation}
The exponential rate of the expected number of labelings of $\Gbconf_{N,d}$ with edge statistics close to a given $\pi$ can be seen to be equal to 
\[ \Sigma(\pi) \defeq \frac{d}{2} H(\pi) - (d-1) H(\pivtx) .\]
For a rigorous proof see \cite[Lemma 4.1 and the proof of Theorem 4]{backhausz2018largegirth}.

Therefore, the first moment method in this context can be phrased as follows.
\begin{claim*}
If $\Sigma(\pi)<0$, then there exists $\eps>0$ such that $\Gb_{N,d}$ a.a.s. has no $\Xc$-labeling with edge statistics $\eps$-close to $\pi$.
\end{claim*}

Furthermore, a second moment argument (supplemented by a standard Azuma-type concentration result) gives the following. 
(Note that Lemma~\ref{lem:boosted_smm} of the next section contains a stronger result that produces labelings satisfying a certain spatial Markov property.)
\begin{claim*} Let $\Pcedge^\Xc$ denote the set of exchangeable distributions on $\Xc \times \Xc$. We say that $\pi \in \Pcedge^\Xc$ satisfies the \emph{second moment condition} if 
\begin{equation} \label{eq:smm_condition}
\Sigma(\nu) \leq 2 \Sigma(\pi) 
\text{ for every coupling $\nu \in \Pcedge^{\Xc \times \Xc}$ of $\pi$ and $\pi$.}
\end{equation}
If $\pi \in \Pcedge^\Xc$ satisfies the second moment condition \eqref{eq:smm_condition}, then, for any $\eps>0$, $\Gb_{N,d}$ a.a.s. has an $\Xc$-labeling with edge statistics $\eps$-close to $\pi$.
\end{claim*}

Many fundamental problems in graph theory can be phrased as questions about whether a graph has a labeling with certain edge statistics. For instance, an independent set of density $\al$ corresponds to a $\{0,1\}$-labeling with edge statistics 
\[ p_{1,1}=0 \quad \text{and} \quad 
p_{0,1}=p_{1,0} = \al .\]
As we will see in Section~\ref{sec:sm_indep}, in this special case we get 
\[ \Sigma(\pi)=\varphi_d(\al)
\quad \text{and} \quad 
\max \Sigma(\nu) 
= \max f_{d,\al} 
=\varphi_d^{(2)}(\al) .\]
%

\section{Typical processes} \label{sec:typical}

When studying labelings or colorings of random regular graphs, a useful tool is to study processes on the limiting graph, which is the $d$-regular infinite tree $T_d$. It is well known that $\Gb_{N,d}$, as well as $\Gbconf_{N,d}$, converges to $T_d$ locally (Benjamini--Schramm convergence) as $N \to \infty$. This is due to the fact that $\Gb_{N,d}$ is essentially large-girth with high probability: for every $\ell \in \Nb$ the number of cycles of length at most $\ell$ is $o(N)$. 

\subsection{Processes}
By an invariant process on $T_d$ we mean a random coloring/labeling, the distribution of which is invariant under the automorphisms of $T_d$. The so-called \emph{tree-indexed Markov chains} form a special class of invariant processes. 
\begin{definition} \label{def:Markov}
Suppose that $\Xc$ is a finite set of colors/labels. A mapping $V(G) \to \Xc$ is a \emph{vertex labeling} of the graph $G$. A random labeling $X \colon V(T_d) \to \Xc$ is a \emph{$T_d$-indexed Markov chain}\footnote{Similar properties are often referred to as \emph{spatial Markov properties}, while the authors of \cite{backhausz2022typicality} call such random labelings \emph{vertex-Markov processes}.} if, conditioned on the label $X_v$ of a vertex $v \in V(T_d)$, $X$ is conditionally independent on the $d$ subtrees obtained after deleting $v$. 
\end{definition}
The distribution $\mu$ of a $T_d$-indexed Markov chain is uniquely described by its marginal $\muedge$ on two neighboring vertices of $T_d$. Note that $\muedge \in \Pcedge^\Xc$, where we write $\Pcedge^\Xc$ for the set of distributions on $\Xc \times \Xc$ that are invariant under exchanging the two coordinates.

\begin{definition} \label{def:typical}
An invariant process $\mu$ is said to be \emph{typical}\footnote{In \cite[Def 1.1]{backhausz2022typicality} the term \emph{strongly typical} was used.} if it can be mimicked on most $d$-regular graphs in the following sense. For every radius $R \in \Nb$ and every positive $\eps$, there exists $N_0=N_0(\eps,R)$ such that for each $N>N_0$ the following holds for $\Gb_{N,d}$ with probability at least $1-\eps$: the graph $\Gb_{N,d}$ has an $\Xc$-labeling with $R$-neighborhood statistics $\eps$-close\footnote{When we talk about the distance of discrete distributions, let us use, say, the total variation distance.} to that of $\mu$.
\end{definition}

\subsection{Entropy}
Given an exchangeable distribution $\pi \in \Pcedge^\Xc$ on $\Xc \times \Xc$, we write $\pivtx$ for its marginal on any of the two coordinates, and define 
\begin{equation} \label{eq:pi_entropy}
\Sigma(\pi) \defeq 
\frac{d}{2} H(\pi) - (d-1) H(\pivtx) .
\end{equation}
See \eqref{eq:Sh_entr} for how Shannon entropy $H(\cdot)$ is defined. The motivation for this definition is that $\Sigma(\pi)$ gives the exponential rate of the expected number of $\Xc$-labelings of $\Gbconf_{N,d}$ with edge statistics close to $\pi$ (see Section~\ref{sec:labelings}). In particular, if $\Sigma(\pi)<0$, then $\Gbconf_{N,d}$ has such a labeling only with exponentially small probability (in $N$). This fact is what we usually refer to as the \emph{first moment method} in random regular graphs. 

When we wish to compute (the exponential rate of) the second moment of the number of such labelings, then we need to determine the expected number of pairs of labelings (one red and one blue), and hence consider $\Sigma(\nu)$ for distributions $\nu$ over
\[ \Xc^4 = 
\big( \textcolor{red}{\Xc} \times 
\textcolor{blue}{\Xc} \big) \times 
\big( \textcolor{red}{\Xc} \times 
\textcolor{blue}{\Xc}\big) .\] 
We will use the following strong form of the \emph{second moment method}.
\begin{lemma}[Backhausz--Bordenave--Szegedy, Thm 1.15 of \cite{backhausz2022typicality}] \label{lem:boosted_smm}
Suppose that $\mu$ is a tree-indexed Markov chain on $T_d$ with the property that for any coupling $\nu \in \Pcedge^{\Xc \times \Xc}$ of $\textcolor{red}{\muedge}$ and $\textcolor{blue}{\muedge}$ we have 
\begin{equation} \label{eq:sm_condition}
\Sigma(\nu) \leq 2 \Sigma(\muedge) .
\end{equation}
Then $\mu$ is (strongly) typical.
\end{lemma}
\begin{remark}
We have to be a bit careful how we think of $\Pcedge^{\Xc \times \Xc}$. Given labels $i,j,k,\ell \in \Xc$, they have the following roles when we write $(\textcolor{red}{i},\textcolor{blue}{j},\textcolor{red}{k},\textcolor{blue}{\ell}) \in \Xc^4$:
\begin{align*}
\textcolor{red}{i} &: 
\quad \text{\textcolor{red}{red} label of the \textbf{first} vertex;}\\
\textcolor{blue}{j} &: 
\quad \text{\textcolor{blue}{blue} label of the \textbf{first} vertex;}\\
\textcolor{red}{k} &: 
\quad \text{\textcolor{red}{red} label of the \textbf{second} vertex;}\\
\textcolor{blue}{\ell} &: 
\quad \text{\textcolor{blue}{blue} label of the \textbf{second} vertex.}
\end{align*}
We need to consider distributions $\nu$ over $\Xc^4$ with the following two properties: 
\begin{itemize}
\item[(i)] $\nu$ is invariant under exchanging the first and second vertex;\footnote{Note that $\nu$ does not have to be invariant under exchanging the two colors.} 
\item[(ii)] the marginal of $\nu$ corresponding to any of the two colors is $\muedge$.
\end{itemize}
We will use the notation 
\[ p_{ij,k\ell} \defeq 
\nu\big( \{ (i,j,k,\ell) \} \big) .\]
Note that the comma separates the labels of the first and second vertex of the edge. With this notation, property (i) can be expressed as  
\begin{equation} \label{eq:nu1}
p_{ij,k\ell} = p_{k\ell,ij} 
\quad (\forall i,j,k,\ell \in \Xc) .
\end{equation}
As for property (ii) regarding the marginals, it translates to 
\begin{align}
\sum_{j,\ell \in \Xc} p_{ij,k\ell} &= p_{i,k} 
\quad (\forall i,k \in \Xc); \label{eq:nu2} \\
\sum_{i,k \in \Xc} p_{ij,k\ell} &= p_{j,\ell} 
\quad (\forall j,\ell \in \Xc). \label{eq:nu3}
\end{align}
\end{remark}

\begin{remark} \label{rem:indep_coupling}
Consider the independent coupling $\nustar \defeq \muedge \times \muedge$ for which $p_{ij,k\ell}=p_{i,k}p_{j,\ell}$ (that is, the red and blue labels are independent). This clearly satisfies the above properties and we have 
\[ \Sigma(\nustar) = 2 \Sigma(\muedge) .\]
Therefore, the second moment condition \eqref{eq:sm_condition} can also be written as $\Sigma(\nu) \leq \Sigma(\nustar), \forall \nu$.

Note that if we choose $\nu$ to be the distribution of two identical copies of $\muedge$ (that is, the blue and red labels are always the same), then $\Sigma(\nu)=\Sigma(\muedge)$. Therefore, if the second moment condition \eqref{eq:sm_condition} holds, then we necessarily have the ``first moment condition'' $\Sigma(\muedge) \geq 0$. 
\end{remark}
%

%
%
%

\subsection{Factors} \label{sec:factors}
As we mentioned in the introduction, many previous lower bounds on $\alstar_d$ used \emph{(randomized) local algorithms} to construct independent sets. These algorithms can be applied to any graphs of essentially large girth, not only to random regular graphs. In the language of graph limits, such local algorithms are equivalent to the so-called \emph{FIID (factor of IID) processes}. See \cite{backhausz2018largegirth} for the precise definition and some basic facts about these processes.

It is easy to see that FIID processes are (strongly) typical. In fact, one can start from a typical process and apply a factor mapping; the resulting process is always typical. More specifically, suppose that $X \colon V(T_d) \to \Xc$ is an invariant process on $T_d$, and let $Z \colon V(T_d) \to [0,1]$ be an IID process (also independent from $X$). In this context, by a factor mapping we mean an $\Aut(T_d)$-equivariant measurable mapping $\Psi \colon \big( \Xc \times [0,1] \big)^{V(T_d)} \to \Yc^{V(T_d)}$. Clearly, $Y \defeq \Psi(X,Z)\colon V(T_d) \to \Yc$ is also invariant under $\Aut(T_d)$. Moreover, it is an easy fact that if $X$ is typical, then so is $Y$.

Essentially, it suffices to consider \emph{block factors}, where the new label of a vertex $v$ is determined by the original labels on $B_R(v)$. In essence, $Y_v$ is computed based on the original labels $X_u$ and the IID labels $Z_u$ for $u \in B_R(v)$, and each vertex has to use the same measurable ``rule''.

In our setting, we will start from some typical Markovian labeling $X_v$, and obtain $Y_v$ using a fixed local algorithm based on the $X$-labels around $v$, possibly using some independent random choices along the way. The resulting labeling is, therefore, also typical.

\section{Second moment for independent sets} \label{sec:sm_indep}

In this section we apply the machinery presented in Section~\ref{sec:typical} to independent sets. 
\begin{definition} \label{def:Markov_indep}
We define an \emph{invariant independent set} as an invariant process on $T_d$ such that each label is $0$ or $1$ with no adjacent $1$'s. If, in addition, the process has the Markov property (see Definition \ref{def:Markov}), then we call it a \emph{Markovian independent set}.

By its \emph{density} we mean the probability that a given vertex has label $1$. Note that a Markovian independent set is uniquely determined by its density.
\end{definition}

Fix some degree $d \geq 3$ and let $\mu = \mu(d,\al)$ denote the distribution of the Markovian independent set on $T_d$ with density $\al$. Our goal in this section is to prove that $\mu$ is typical provided that $d$ and $\al$ satisfy a certain second moment condition.

The vertex and edge marginals of $\mu$ are as follows:
\begin{align*} 
\muvtx&: \quad 
p_1=\al; \; p_0=1-\al; \\
\muedge&: \quad 
p_{1,1}=0; \; p_{0,1}=p_{1,0}=\al; \; p_{0,0}=1-2\al .
\end{align*}
Now let $\nu$ be an arbitrary coupling of $\muedge$ and $\muedge$. (One may think of two random independent sets, one red and one blue, both of density $\al$.) As explained in Section~\ref{sec:typical}, we need to consider distributions $\nu$ over $\{0,1\}^4$ for which the probabilities $p_{ij,k\ell}$ satisfy \eqref{eq:nu1}, \eqref{eq:nu2}, \eqref{eq:nu3}.

Let $\be \defeq p_{11}$. (This corresponds to the density of the intersection of the blue and red independent sets. Note that this is different from $p_{1,1}=0$, which is the probability that two neighboring vertices both have $0$ blue labels.) Then we have 
\[ \nuvtx: \quad 
p_{11}=\be; \; p_{01}=p_{10}=\al-\be; \; p_{00}=1-2\al+\be .\]
It follows that 
\[ p_{00,11}=p_{11,00}=\be .\]
Next let $\ga \defeq p_{10,01}=p_{01,10}$. Then 
\[ p_{01,00}+p_{01,10}=p_{01}=\al-\be 
\text{, and hence } 
p_{01,00}= \al-\be-\ga = p_{00,01}.\]
Similarly,
\[ p_{10,00}+p_{10,01}=p_{10}=\al-\be 
\text{, and hence } 
p_{10,00}= \al-\be-\ga = p_{00,10}.\]
Then the remaining probability is $p_{00,00}=1-4\al+2\be+2\ga$, and we get that 
\begin{align*}
H(\nuvtx) &=  h(\be)+ 2h(\al-\be) + h(1-2\al+\be) ; \\
H(\nu) &= 2h(\be)+2h(\ga)+4h(\al-\be-\ga) + h(1-4\al+2\be+2\ga) .
\end{align*} 
In order to check the second moment condition of Lemma~\ref{lem:boosted_smm}, we need to maximize 
\[ \Sigma(\nu) = \frac{d}{2} H(\nu) - (d-1) H(\nuvtx) \]
and show that the maximum is attained at the independent coupling $\nustar$ (see Remark~\ref{rem:indep_coupling}). For fixed $d$ and $\al$, this expression is a two-variable function in $\be$ and $\ga$ that we denote by $\fhat_{d,\al}$:
\begin{align}
\begin{split} \label{eq:fhat}
\fhat_{d,\al}(\be,\ga)  
&\defeq \frac{d}{2} \bigg( 2h(\be)+2h(\ga)
+4h(\al-\be-\ga) + h(1-4\al+2\be+2\ga) \bigg) \\
& - (d-1) \bigg( h(\be)+ 2h(\al-\be) 
+ h(1-2\al+\be) \bigg) .
\end{split}
\end{align} 
Note that second line does not depend on $\ga$, while the first line is concave in $\ga$ because $h$ is a concave function. Therefore, for any fixed $\be$, $\fhat_{d,\al}(\be,\ga)$ is concave in $\ga$, and takes its maximum, where the partial derivative w.r.t. the second variable vanishes. We have  
\[ \partial_2 \fhat_{d,\al} (\be,\ga)
= d \bigg( -\log \ga +2 \log(\al-\be-\ga) - \log(1-4\al+2\be+2\ga) \bigg) ,\]
which vanishes if and only if 
\[ (\al-\be-\ga)^2=\ga(1-4\al+2\be+2\ga) .\]
This is a quadratic equation in $\ga$ with one positive root that we denote by 
\begin{align}
\Ga(\be) 
&\defeq -(1/2-\al)+\sqrt{(1/2-\al)^2+(\al-\be)^2} 
\label{eq:Ga_first} \\
&= \frac{(\al-\be)^2}{(1/2-\al) 
+\sqrt{(1/2-\al)^2+(\al-\be)^2}} .
\label{eq:Ga_second}
\end{align}  
In conclusion, $\fhat_{d,\al}(\be,\ga)$ is maximized in the second variable at $\ga=\Ga(\be)$. We define the one-variable function $f_{d,\al}$ as follows:
\begin{equation} \label{eq:f_al}
f_{d,\al}(\be) 
\defeq \max_{\ga} \fhat_{d,\al}(\be,\ga) 
= \fhat_{d,\al}\big(\be,\Ga(\be) \big) .
\end{equation}

There are two important particular substitutions, namely, $\be=\al$, $\ga=0$ (identical coupling) and $\be=\ga=\al^2$ (independent coupling). Next we show that, in these special cases, the value of $\fhat_{d,\al}$ can be expressed using the first moment rate-function $\varphi_d(\al)$ defined in \eqref{eq:phi}. 
\begin{claim} \label{cl:special}
We have 
\begin{align*}
\Ga(\al) =0 \text{ and }
f_{d,\al}(\al) &= 
\fhat_{d,\al}(\al,0) = \varphi_d(\al) ; \\
\Ga(\al^2) = \al^2 \text{ and }
f_{d,\al}(\al^2) &= 
\fhat_{d,\al}(\al^2,\al^2) = 2 \varphi_d(\al) .
\end{align*}
\end{claim}
\begin{proof}
Due to \eqref{eq:Ga_second} we clearly have $\Ga(\al)=0$, while $\Ga(\al^2)=\al^2$ follows from \eqref{eq:Ga_first} and the identity 
\begin{equation} \label{eq:identity}
\big( 1/2 - \al \big)^2 
+ \big( \al - \al^2 \big)^2 
= \big( 1/2 - \al + \al^2 \big)^2 .
\end{equation}
Then $f_{d,\al}(\al)=\varphi_d(\al)$ and $f_{d,\al}(\al^2)=2\varphi_d(\al)$ follow easily from \eqref{eq:fhat} if we make the particular substitutions.
\end{proof}

Our aim is to apply Lemma~\ref{lem:boosted_smm} in the above setting. Here the independent coupling corresponds to the probabilities 
\[ \be = 
p_1^2=\al^2 
\quad \text{and} \quad 
\ga = 
p_{1,0}p_{0,1} = \al^2 .\]
Therefore, the condition of the lemma is that $\fhat_{d,\al}$ is maximized at $\be=\ga=\al^2$. According to Claim~\ref{cl:special}, this is the same as saying that $f_{d,\al}$ is maximized at $\be=\al^2$, and the following result is obtained. 
\begin{theorem} \label{thm:Markovian_indep}
Let $d \geq 3$ be an integer and let $\al \in (0,1/2)$. If condition \eqref{eq:the_cond} holds (that is, $f_{d,\al}$ attains its maximum at $\al^2$), then $\mu(d,\al)$, the Markovian independent set of density $\al$, is a typical process on $T_d$.

Note that, using Definition~\ref{def:alsm}, we can simply state this result as follows: 
\begin{equation*}
\mu(d,\alsm_d) 
\text{ is typical; in particular, }
\alstar_d \geq \alsm_d .
\end{equation*}
\end{theorem}

\begin{figure}
\centering
\includegraphics[width=0.5\linewidth]{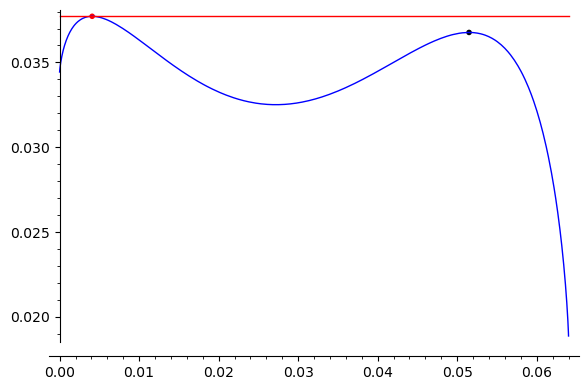}
\caption{The plot of $f_{d,\al}(\be)$ for $d=100$ and $\al=0.064$. The red dot corresponds to $\be=\al^2$ (independent coupling). The second moment method works here because the red dot is above the second local maximum (black dot) so it is the global maximum.}
\label{fig:f_plot}
\end{figure}

The question arises: what does the graph of $f_{d,\al}$ look like, and how can we check if condition \eqref{eq:the_cond} holds for specific $d$ and $\al$? In Figure~\ref{fig:f_plot} we plotted $f_{d,\al}(\be)$ on $[0,\al]$ in the case $d=100$ and $\al=0.064$. The first local maximum corresponds to $\be=\al^2$. Since it is larger than the second local maximum, we can conclude that \eqref{eq:the_cond} is satisfied, and hence the Markovian independent set of density $0.064$ is typical on $T_{100}$.

Do we always see a similar picture? How does the plot of $f_{d,\al}$ change as $\al$ varies? We encounter the following behavior\footnote{We will not prove all of the claimed behavior; only the pieces necessary for proving Theorem~\ref{thm:vs_alfm} are included in the paper. See Section~\ref{sec:local_maxima} for precise statements and rigorous proofs.} for each $d \geq 3$. There are two density thresholds (the second of which can be expressed with an explicit formula):
\[ 0<\al_d^{\textrm{1st}} 
< \al_d^{\textrm{2nd}}
=\frac{\sqrt{d-1}-1}{d-2} 
< \frac12 .\]
There are three phases:
\begin{itemize}
\item $0< \al < \al_d^{\textrm{1st}}$:\\
$f_{d,\al}$ has a single local maximum at $\al^2$, which is, therefore, also the global maximum, and hence the second moment condition \eqref{eq:the_cond} is \textbf{always satisfied}.
\item $\al_d^{\textrm{1st}}<\al<\al_d^{\textrm{2nd}}$:\\
$f_{d,\al}$ has two local maxima: one at $\be_1^{\max}=\al^2$ and one at some $\be_2^{\max}>\al^2$. The second moment condition \eqref{eq:the_cond} is \textbf{satisfied if and only if} 
\begin{equation} \label{eq:first_is_larger}
f_{d,\al}\big( \be_1^{\max} \big) 
\geq f_{d,\al}\big( \be_2^{\max} \big) .
\end{equation}
\item $\al_d^{\textrm{2nd}} < \al < 1/2$:\\
$f_{d,\al}$ has a local \underline{minimum} at $\be=\al^2$, so the second moment condition \eqref{eq:the_cond} is \textbf{never satisfied}.
\end{itemize}

Therefore, in the intermediate phase (between $\al_d^{\textrm{1st}}$ and $\al_d^{\textrm{2nd}}$), checking condition \eqref{eq:the_cond} amounts to finding $\be_2^{\max}$ and verifying \eqref{eq:first_is_larger}. (For specific $d$ and $\al$, this can be done numerically with a computer very quickly.) 

We finish this section with the following simple observation.
\begin{remark}
Due to Claim~\ref{cl:special} we have  
\[ \eqref{eq:the_cond} 
\Rightarrow 
\underbrace{f_{d,\al}(\al^2)}_{=2\varphi_d(\al)} 
\geq 
\underbrace{f_{d,\al}(\al)}_{=\varphi_d(\al)} 
\Leftrightarrow \varphi_d(\al) \geq 0 
\Leftrightarrow \al \leq \alfm_d . \]
Therefore, $\alsm_d \leq \alfm_d$ (which also follows from the fact that $\alstar_d$ is sandwiched between the two). 
\end{remark}
%

\section{Augmented independent sets} \label{sec:augment_indep}

In this section we discuss the local improvements one can make to an independent set with the spatial Markov property. Let us start from a Markovian independent set of density $\al$ (see Definition~\ref{def:Markov_indep}). If we have two neighboring vertices, then the conditional probability that one has label $0$ conditioned on the other having label $0$ is equal to 
\[ p \defeq \frac{1-2\al}{1-\al} .\]
We call a vertex \emph{full-zero} if its label is $0$ and the labels of all its neighbors are also $0$. Then, for any vertex $v$ of $T_d$, we clearly have
\begin{equation} \label{eq:prob_fullzero}
\Pb( v \text{ is full-zero} ) 
= (1-\al) p^d .
\end{equation}
We will also need the following conditional probability. If $u$ and $v$ are adjacent vertices, then 
\begin{equation} \label{eq:cond_full_zero}
\Pb( u \text{ is full-zero} \, | \, v \text{ is full-zero} ) 
= p^{d-1} .
\end{equation}
It follows that the probability of a vertex being an isolated full-zero vertex (that is, a full-zero vertex with no full-zero neighbor) is
\begin{equation} \label{eq:prob_isolated}
\Pb( v \text{ is full-zero and no neighbor of $v$ is full-zero} )
= (1-\al) p^d \big( 1- p^{d-1} \big)^d .
\end{equation}
Now consider the induced subgraph over all full-zero vertices. We will need that the connected components of this subgraph are finite (almost surely), that is, the corresponding Galton--Watson process dies out almost surely. The latter happens if and only if the conditional probability \eqref{eq:cond_full_zero} is at most $1/(d-1)$. Next we check that this condition holds (unless $\al$ is too small) so that we may safely assume that the components are finite.
\begin{claim} \label{cl:p_bnd}
\[ \text{If } 
\al \geq \frac{\log(d-1)}{d-1}
\text{, then }
p^{d-1} < \frac{1}{d-1} .\]
\end{claim}
\begin{proof}
Note that 
\[ p = \frac{1-2\al}{1-\al} < 1-\al \leq 1 - \frac{\log(d-1)}{d-1} 
\leq \exp\left( - \frac{\log(d-1)}{d-1} \right).\]
Therefore, 
\[ p^{d-1} < \exp\left( - \log(d-1) \right) = \frac{1}{d-1} .\]    
\end{proof}
It is easy to design a randomized local algorithm that selects an independent set inside each finite component with size at least half of the component.\footnote{The components are trees and hence bipartite. Take a ``canonical'' proper 2-coloring and choose the larger color class. If they have the same size, choose the one in which the sum of IID labels is larger.} Let us add these to the original independent set. Of course, the resulting set remains to be independent. What can we say about its density? We augmented the original set with all the one-vertex components (i.e., isolated full-zero vertices) and at least half of the remaining full-zero vertices. Therefore, using \eqref{eq:prob_fullzero} and \eqref{eq:prob_isolated}, the new density is at least 
\begin{equation} \label{eq:alh}
\alh \defeq \al + 
(1-\al)p^d\frac{1+\big( 1- p^{d-1} \big)^d}{2} .
\end{equation}
We summarize our findings in the next theorem.
\begin{theorem} \label{thm:augmented_indep}
Let $d \geq 3$ be an integer and $\al \in (0,1/2)$, and set $p=(1-2\al)/(1-\al)$. Consider the function $f_{d,\al} \colon [0,\al] \to \Rb$ defined in \eqref{eq:f_al} and suppose that it attains its maximum at $\al^2$. If $(d-1)p^{d-1} \leq 1$, then there exists a typical independent set on $T_d$ with density $\alh$ as in \eqref{eq:alh}.
\end{theorem}
Setting $\al=\alsm_d$, we get the boosted bound $\alh=\alaug_d$ stated in \eqref{eq:al_aug}. It can be shown that the condition $(d-1)p^{d-1} \leq 1$ is always satisfied\footnote{Why can we assume that the condition $(\ast)\; (d-1)p^{d-1} \leq 1$ (for $\al=\alsm_d$) holds true? In this paper we only use Corollary~\ref{cor:alaug} for producing numerical bounds, in which case it suffices to check $(\ast)$ for the given $d$ (which we did). To show that $(\ast)$ holds true for every $d$, one may combine Theorem~\ref{thm:vs_alfm} and Theorem~\ref{thm:alfm_bnd} to prove that $\alsm_d \geq \frac{\log(d-1)}{d-1}$ for sufficiently large $d$. This, in turn, implies $(\ast)$ by Claim~\ref{cl:p_bnd}.} for this $\al$, and hence we get the following corollary. 
\begin{corollary} \label{cor:alaug}
For $d \geq 3$ let $\alaug_d$ be the density defined in \eqref{eq:al_aug}. Then $T_d$ has a typical independent set of density $\alaug_d$. In particular, $\alstar_d \geq \alaug_d$.
\end{corollary}
Table~\ref{tab:bound_comparision} in the introduction compares the original density $\alsm_d$ and the increased density $\alaug_d$ for specific degrees.

One could achieve a slightly higher density by also computing the probability of components with three vertices (two of which could be added to the independent set), or by considering further local improvements, but these would only lead to marginal gains.

\subsection{The asymptotic gain of the boost} \label{sec:asymptotic_gain}
The above local improvement increases the density substantially for small degrees. Asymptotically, however, the gain is only of order $1/d^{2-o(1)}$. More precisely, one can show that 
\begin{equation} \label{eq:gain}
\alaug_d < \alsm_d + \frac{(\log d)^2}{d^2} 
= \alsm_d + \frac{1}{d^{2-o(1)}} .
\end{equation}
From \eqref{eq:alstar_vs_alfm} and \eqref{eq:alsm_vs_alfm} we have 
\[ \alsm_d = \alfm_d - \frac{1}{d^{3/2-o(1)}} 
\quad \text{and} \quad 
\alstar_d = \alfm_d - \frac{1}{d^{2-o(1)}} .\]
It means that (for large $d$) we do not get significantly closer to the true value with the augmented bound $\alaug_d$.

We finish this section by sketching an argument for proving \eqref{eq:gain}. As we pointed out in the proof of Claim~\ref{cl:p_bnd}, we have $p < 1-\al$. Moreover, the probability \eqref{eq:prob_fullzero} that a vertex is full-zero is clearly an upper bound for the density gain:
\[ \hat{\al}-\al \leq (1-\al)p^d \leq p^d 
< (1-\al)^d \leq \exp(-d\al) .\]
From Theorem~\ref{thm:vs_alfm} and Theorem~\ref{thm:alfm_bnd}(c), one can easily deduce that 
\[ \alsm_d \geq \frac{2 \log d}{d} - \frac{2 \log\log d}{d} .\]
Therefore, for $\al=\alsm_d$ the gain is less than  
\[ \exp\big( -2 \log d + 2 \log\log d \big) 
= \frac{(\log d)^2}{d^2} .\]
%

\section{Star decompositions} \label{sec:star_decomp}

To further demonstrate the strength of our approach, we consider the following problem. Given an integer $d/2 < k < d$, can the edges of $\Gb_{N,d}$ be partitioned into edge-disjoint stars, each containing $k$ edges? We restrict ourselves to those $N$ for which $Nd/2$ (i.e., the number of edges) is divisible by $k$. If such a partition exists with probability $1-o_N(1)$, then we say that $\Gb_{N,d}$ asymptotically almost surely (a.a.s.) has a $k$-star decomposition.


The case $d=4, k=3$ was settled (affirmatively) in \cite{delcourt2018random}. Later, that method was extended in \cite{delcourt2025decomposing}, and a sufficient condition was given for general pairs $d,k$. Both papers used a second moment method and checking the required conditions led to very technical computations.

The problem is closely related to the asymptotic independence ratio $\alstar_d$ of $\Gb_{N,d}$. More specifically, if $\Gb_{N,d}$ a.a.s.~has a $k$-star decomposition, then 
\[ \alstar_d \geq 1 - \frac{d}{2k} 
\text{, or equivalently: } 
k \leq \frac{d}{2(1-\alstar_d)} .
\]
This gives a necessary condition for the existence of a $k$-star decomposition. It was conjectured in \cite{delcourt2025decomposing} that this is essentially the only requirement. This was confirmed for most degrees (in an asymptotic sense) in \cite{harangi2025star1}. In that paper we showed that the existence of large independent sets with an additional property (called \emph{thinness}) actually implies the existence of star decompositions provided that the graph has some properties regarding the edge density of its induced subgraphs. We will use this general result here. First we define the notion of thin sets in graphs.
\begin{definition} \label{def:thin}
Given a positive integer $\dhat < d$, we say that a set $A \subset V(G)$ is \emph{$\dhat$-thin in $G$} if every outside vertex $v \in V(G) \sm A$ has at most $\dhat$ neighbors in $A$.
\end{definition}
Then we have the following sufficient condition for the existence of a $k$-star decomposition.
\begin{lemma}[\cite{harangi2025star1}] \label{lem:sd1}
Let $G$ be a $d$-regular graph on $N$ vertices and let $k>d/2$. Set 
\[ \al_{d,k} \defeq 1 - \frac{d}{2k} .\]
For $U \subseteq V(G)$ let $e[U]$ denote the number of edges of $G$ inside $U$. Suppose that $A \subset V(G)$ satisfies the following conditions for a positive integer $\dhat < k$ and a real number $0<c<1$.
\begin{enumerate}[(i)]
\item $A$ is an independent set of density $\al_{d,k}$.
\item $A$ is $\dhat$-thin.
\item For any set $U \subseteq A^\comp = V(G) \sm A$ with $|U|\leq cN$ we have
\[ e[U] \leq (d-k) |U| .\]
\item For any set $W \subseteq A^\comp$ with $|W| < (1-\al_{d,k}-c)N$ we have 
\[ e[W] \leq (k-\dhat) |W| .\]
\end{enumerate}
Then $G$ has a $k$-star decomposition.
\end{lemma}
If we choose the parameters properly, then conditions (iii) and (iv) actually hold true a.a.s.~for random regular graphs. In fact, one can simply check the numerical condition of \cite[Cor 3.4]{harangi2025star1} to check (iii) and (iv) in specific cases. We omit the details here as this is not the main focus of the current paper. Instead, we focus on the first two conditions.

To find a set $A$ satisfying (i) and (ii), we start from a Markovian independent set of density $\al$. As before, $p \defeq (1-2\al)/(1-\al)$ is the conditional probability that a vertex has label $0$ conditioned on the event that a fixed neighbor has label $0$. Then, by the Markov property, the probability that a vertex with label $0$ has $\ell$ neighbors with label $1$ is the following: 
\[ \binom{d}{\ell} (1-p)^\ell p^{d-\ell} .\]
For any vertex $v \notin A$ with $\ell > \dhat$ neighbors in $A$, we may remove $\ell-\dhat$ vertices from $A$ to ensure that $v$ does not violate the thinness condition. (Note that, although we created some new outside vertices during this removal process, they all have $0$ neighbors in $A$.) This simple local algorithm removes a subset of $A$ of density at most 
\[ (1-\al) \sum_{\ell=\dhat+1}^d (\ell-\dhat) 
\binom{d}{\ell} (1-p)^\ell p^{d-\ell} .\]
Therefore, we have proved the following.
\begin{lemma} \label{lem:thinning}
Assume that the Markovian independent set of density $\al$ is typical, and let $p \defeq (1-2\al)/(1-\al)$. For $\dhat < d$ there exists a $\dhat$-thin invariant independent set $A_{\mathrm{thin}} \subseteq A$ of density at least
\[ \althin \defeq 
\al - (1-\al) \sum_{\ell=\dhat+1}^d (\ell-\dhat) 
\binom{d}{\ell} (1-p)^\ell p^{d-\ell} ,\]
which is also typical.
\end{lemma}

For a specific example, consider the degree $d=33$. In this case, the density of the Markovian independent set is $\al \approx 0.13195395$, while the corresponding conditional probability is $p \approx 0.84798738$. Then, by Lemma~\ref{lem:thinning} we get the following densities for $\dhat$-thin independent sets for $\dhat=8,9,\ldots,14$:

\medskip

\begin{center}
\begin{tabular}{r||l|l|l}
$\dhat$ & deleted & $\althin$ & star decomposition \\
\hline
$8$  & $0.073901$ & $0.058052$ & implies $k=17$ with $c=0.1$  \\
$9$  & $0.027726$ & $0.104227$ & implies $k=18$ with $c=0.1$   \\
$10$ & $0.009286$ & $0.122667$ &    \\
$11$ & $0.002778$ & $0.129175$ &    \\
$12$ & $0.000743$ & $0.131210$ &    \\
$13$ & $0.000177$ & $0.131776$ & implies $k=19$ with $c=0.1$ \\
$14$ & $0.000038$ & $0.131915$ &  
\end{tabular}
\end{center}

\medskip
Then one can use Lemma~\ref{lem:sd1} to prove the existence of $k$-star decompositions in $\Gb_{N,33}$ for $k=17,18,19$. (Note that star decompositions cannot exist for $k \geq 20$ because $\alstar_{33} \leq 0.14401 < 0.175 = 1- \frac{33}{40}$.) As we mentioned, conditions (iii) and (iv) can be easily checked (for the choice $c=0.1$) with the tools in \cite{harangi2025star1}. This way one can reproduce many results proven in \cite{delcourt2025decomposing} with substantially smaller technical difficulty.

\section{Estimates} \label{sec:estimates}
In this section we collected some general estimates that we will need later on. We start with $\log$'s power series around $1$:
\begin{equation} \label{eq:log_power_series}
\log(1-x) = - \sum_{k=1}^\infty \frac{x^k}{k} 
\quad \text{for } x \in [0,1) .
\end{equation} 
In particular, we have the following bounds: 
\begin{align} 
\log(1-x) &\leq -x - \frac12 x^2 \leq -x 
\quad \text{for } x \in [0,1) ; 
\label{eq:log_upp_bnd}\\
\log(1-x) &\geq -x - x^2 
\quad \text{for } x \in [0,2/3] ; 
\label{eq:log_low_bnd1}
\end{align}

Now we turn to the function $h(1-x)=-(1-x)\log(1-x)$, which can be expressed as follows:
\begin{equation} \label{eq:h_ps}
h(1-x) = x - \sum_{k=2}^\infty \frac{x^k}{k(k-1)} 
= x - \frac{x^2}{2} - \frac{x^3}{6} - \frac{x^4}{12} - \cdots \quad \text{for any } x \in [0,1) .
\end{equation}
It follows that 
\begin{equation} \label{eq:h_diff2_simple}
h(1-x-\eps) - h(1-x) 
= \eps - x\eps -\frac12 \eps^2
- \sum_{k=3}^\infty \frac{(x+\eps)^k-x^k}{k(k-1)} .
\end{equation}
Furthermore, the derivative 
\[ h'(x) = -1 - \log(x) \]
is monotone decreasing in $x$, and hence 
\begin{equation} \label{eq:h_diff1}
h(x-\eps) - h(x) 
= - \int_{x-\eps}^x h'(z) \,\mathrm{d}z 
\leq -\eps h'(x) = \eps(1+\log x) 
\quad \text{for any } 0 \leq \eps \leq x \leq 1.
\end{equation}
Next we consider functions of the form $h(x)+cx$ and $h(x)+cx+Cx$, and study their maximum.
\begin{claim} \label{cl:hx+cx}
Let $c<0$. The function 
\[ f(x) \defeq h(x)+cx = x \big( -\log(x) + c \big) 
\quad (0 \leq x \leq 1) \]
is maximized at $x_0 \defeq \exp(c-1)$ with maximum value $f(x_0)=x_0$. 
\end{claim}
\begin{proof}
The first and second derivative of $f$:
\begin{align*}
f'(x) &= -1 -\log(x) + c ;\\
f''(x) &= -\frac{1}{x} .
\end{align*} 
So $f$ is a concave function attaining its maximum at a point $x_0$ where $f'$ vanishes:  
\[ f'(x_0) = 0 \; \Leftrightarrow 
\log(x_0)=c-1 \; \Leftrightarrow 
x_0= \exp(c-1) .\]
Furthermore, for $x_0$ we have 
\[ f(x_0)= x_0 \big( \underbrace{-\log(x_0) + c}_{=1} \big) 
= x_0 ,\]
and the proof is complete.
\end{proof}

In fact, we will need the following corollary.
\begin{claim} \label{cl:hx+cx+Cxx}
Let $c<0$ and $C>0$ such that 
\[ 2C < \exp(-c) .\]
Consider the function 
\[ g(x) \defeq h(x)+cx+Cx^2 = x \big( -\log(x) + c + Cx \big) .\]
Then 
\[ \max_{0 \leq x \leq 1/(2C)} g(x) 
\leq \exp(c-1) + C\exp(2c) .\]
\end{claim}
\begin{proof}
The first and second derivative of $g$:
\begin{align*}
g'(x) &= -1-\log x + c + 2Cx ;\\
g''(x) &= -\frac{1}{x} +  2C .
\end{align*}
We conclude that $g$ is concave over $(0,\frac{1}{2C}]$. We define 
\[ x_0 \defeq \exp(c-1) 
\quad \text{and} \quad 
x_1 \defeq e x_0 = \exp(c) .\]
Then, using the condition $2C < \exp(-c)$, we get that 
\[ g'(x_1) = -1-c+c+2C \exp(c) = -1 + 2C \exp(c) < 0 .\]
Therefore, the derivative is negative at $x_1$. Note that, due to the same condition, $x_1$ lies in the interval $(0,\frac{1}{2C}]$, where $g$ is concave. Consequently, the maximum of $g$ (over this interval) must be attained at some $0 < x < x_1$. Using Claim~\ref{cl:hx+cx} we conclude that 
\begin{align*}
\max_{0 \leq x \leq 1/(2C)} g(x) 
= \max_{0 \leq x \leq x_1} g(x)
&\leq \bigg( \max_{0 \leq x \leq x_1} h(x)+cx \bigg) 
+ \bigg( \max_{0 \leq x \leq x_1} Cx^2 \bigg) \\
&= x_0 + C x_1^2 = \exp(c-1)+ C\exp(2c) ,
\end{align*}
and the proof is complete.
\end{proof}
%

\section{Local maxima} \label{sec:local_maxima}

In this section we analyze the function $f_{d,\al}(\be)$, mainly focusing on its derivatives and local optima. Recall that $f_{d,\al}(\be) = \fhat_{d,\al}(\be,\Ga(\be))$; see Section~\ref{sec:sm_indep} for the details. 

\subsection{The derivatives} \label{sec:derivatives}
We use $\partial_1$ and $\partial_2$ to denote the partial derivatives w.r.t. the first and second variable, respectively. Then we have the following for the derivative of $f_{d,\al}$:
\[ f'_{d,\al}(\be) 
= \partial_1 \fhat_{d,\al} 
\big( \be, \Ga(\be) \big) 
+ \underbrace{\partial_2 \fhat_{d,\al} 
\big( \be, \Ga(\be) \big)}_{=0} 
\cdot \Ga'(\be) 
= \partial_1 \fhat_{d,\al} 
\big( \be, \Ga(\be) \big)  ,\]
where 
\begin{align*}
\partial_1 \fhat_{d,\al} (\be,\ga)
&= d \bigg( -\log \be + 
\underbrace{2 \log(\al-\be-\ga) 
- \log(1-4\al+2\be+2\ga)}_{(\ast)} \bigg) \\
&- (d-1) \bigg( -\log \be + 2 \log(\al-\be)-\log(1-2\al+\be) \bigg) .
\end{align*}  
If $\ga=\Ga(\be)$, then the expression $(\ast)$ simplifies to $\log \ga = \log \Ga(\be)$, and we have 
\begin{align} 
\begin{split} \label{eq:f_der}
f'_{d,\al}(\be) 
&= \partial_1 \fhat_{d,\al} 
\big( \be, \Ga(\be) \big) \\
&= d \log \Ga(\be) - \log \be 
- (d-1) \bigg( 2 \log(\al-\be)-\log(1-2\al+\be) \bigg) \\
& \stackrel{\eqref{eq:Ga_second}}{=} -d \log\left( 1/2-\al 
+\sqrt{(1/2-\al)^2+(\al-\be)^2} \right) \\
&+ (d-1)\log(1-2\al+\be)
+2\log(\al-\be) 
-\log \be .
\end{split}
\end{align}
It follows easily that  
\begin{align}
\begin{split} 
\lim_{\be \to 0+} f'_{d,\al}(\be) &= +\infty ; \\
\lim_{\be \to \al-} f'_{d,\al}(\be) &= -\infty .
\label{eq:endpoint_der}
\end{split}
\end{align}
In particular, $f_{d,\al}$ is monotone increasing near $0$ and monotone decreasing near $\al$.

As for the second derivative, we have  
\begin{align} 
\begin{split} \label{eq:f_derder}
f''_{d,\al}(\be) 
&= \frac{d(\al-\be)}{G_\al(\be)} 
+ \frac{d-1}{1-2\al+\be}
-\frac{2}{\al-\be} 
-\frac{1}{\be} \text{, where} \\
G_\al(\be) &= \big( 1/2-\al \big) 
\sqrt{(1/2-\al)^2+(\al-\be)^2} 
+(1/2-\al)^2+(\al-\be)^2
\end{split}
\end{align}
Then we have the following for the third derivative:
\begin{equation} \label{eq:f_derderder}
f'''_{d,\al}(\be) = 
- \frac{d}{G_\al(\be)} 
- \frac{d(\al-\be) G'_\al(\be)}{G^2_\al(\be)} 
- \frac{d-1}{(1-2\al+\be)^2}
- \frac{2}{(\al-\be)^2}
+ \frac 1 {\be^2} .
\end{equation}

\subsection{Crude asymptotics} \label{sec:crude}
We will not need this when we do our precise calculations later on; nevertheless let us quickly do some crude asymptotics to get a rough picture of what the graph of $f_{d,\al}$ looks like for large $d$. 

In what follows, we assume that $\al =\Oc\big( \frac{\log d}{d} \big)$ and write $\approx$ when an equality holds with an (additive) error $o(\log)$. Then, using \eqref{eq:f_der}, we have   
\begin{align*}
f'_{d,\al}(\be) &= 
2d\al - d(2\al-\be) + 2\log(\al-\be) - \log \be 
+ \Oc(d\al^2)\\
&\approx d\be + 2\log(\al-\be) - \log \be .
\end{align*}
Furthermore, if $1/d \leq \be \leq \al$, then $\log \be \approx \log d$. Similarly, if $1/d \leq \al-\be \leq \al$, then  $\log (\al-\be) \approx \log d$. We distinguish the following three regimes for $\be$.
\begin{itemize}
\item If $0 \leq \be \leq \frac{1}{d}$, then $d \be = \Oc(1) \approx 0$ and we get  
\[ f'_{d,\al}(\be) \approx 0 -2 \log d - \log \be ,\]
meaning that $f_{d,\al}$ has a local maximum at some %
\[ \be_1^{\max}=\frac{1}{d^{2+o(1)}} .
\quad \text{(In fact, we have } \be_1^{\max}=\al^2 .) 1\]
\item If $\ds \frac{1}{d} \leq \be \leq \al - \frac{1}{d}$, then   
\[ f'_{d,\al}(\be) \approx d \be - 2\log d + \log d ,\]
which vanishes when $d \be \approx \log d$. It follows that $f_{d,\al}$ has a local minimum at some 
\[ \be^{\min} = (1+o(1)) \frac{\log d}{d} \] 
provided that $\al> \frac{\log d}{d}$.
%
%
\item If $\ds \al - \frac{1}{d} \leq \be \leq \al$, then $d \be = d \al + \Oc(1) \approx d \al$, and we obtain 
\[ f'_{d,\al}(\be) \approx d \al + 2\log(\al-\be) + \log d .\]
If $\al = (2+o(1)) \frac{\log d}{d}$, then we have another local maximum at some 
\[ \be_2^{\max} = \al - \frac{1}{d^{3/2+o(1)}} . \]
\end{itemize}

From here we could actually quickly deduce that $\alsm_d = \alfm_d - 1/d^{3/2+o(1)}$. However, if we want to get the precise statement of Theorem~\ref{thm:vs_alfm}, then we will have to perform very delicate estimations. This is what the rest of the paper is devoted to.

\subsection{Function behavior} \label{sec:behavior}
In Section~\ref{sec:sm_indep} we described the behavior of $f_{d,\al}$ in different density regimes. Our next goal is to rigorously prove those aspects of this behavior that will be useful in the proof of Theorem~\ref{thm:vs_alfm}. First we show that $f_{d,\al}$ always has a stationary point at $\al^2$, and it is a local maximum if $\al$ is below the second density threshold $\al_d^{\textrm{2nd}}$. 
\begin{claim} \label{cl:al_sq}
For every $d \geq 3$ and every $\al \in (0,1/2)$, the function $f_{d,\al}$ has a stationary point at $\be=\al^2$. Furthermore, this is a local maximum provided that  
\begin{equation} \label{eq:alsq_cond}
\al < \al_d^{\textrm{2nd}} \defeq \frac{\sqrt{d-1}-1}{d-2} .
\end{equation}
We remark that $2(\log d)/d<\al_d^{\textrm{2nd}}$ for all $d \geq 103$ (and $\alfm_d < \al_d^{\textrm{2nd}}$ for all $d \geq 29$).
\end{claim}
\begin{proof}
Recall that we computed the derivatives $f'_{d,\al}$ and $f''_{d,\al}$ in Section~\ref{sec:derivatives}. Using the identity \eqref{eq:identity} we get that 
\[ 1/2-\al +\sqrt{(1/2-\al)^2+(\al-\be)^2} \text{ at } \be=\al^2 \text{ is simply } 1-2\al+\al^2=(1-\al)^2 ,\] 
and hence 
\begin{align*}
f'_{d,\al}(\al^2) &= -
\underbrace{\log( 1-2\al + \al^2 )}_{=2\log(1-\al)} 
+ 2 \underbrace{\log(\al-\al^2)}_{=\log \al + \log(1-\al)} - \underbrace{\log(\al^2)}_{=2\log \al}
= 0 ;\\
f''_{d,\al}(\al^2) &=  
\frac{d\al(1-\al)}{(1-\al)^2 \big( 1/2 - \al + \al^2 \big)}
+\frac{d-1}{(1-\al)^2}
-\frac{2}{\al(1-\al)}
-\frac{1}{\al^2} \\
&= \frac{(d-2)\al^2 + 2\al -1}
{\al^2(1-\al)^2(1-2\al+2\al^2)}.
\end{align*}
It follows that $\al^2$ is always a stationary point. Moreover, it is a local maximum if 
\[ f''_{d,\al}(\al^2)<0 \Leftrightarrow
(d-2)\al^2 + 2\al -1 < 0 ,\]
which can be checked to be equivalent to \eqref{eq:alsq_cond}.
\end{proof}

We continue with analyzing the local maxima of $f_{d,\al}$. We will assume that the density $\al$ is in the relevant region.
\begin{lemma} \label{lem:local_maxima}
Suppose that $\al \in [\al_d^-,\al_d^+]$, where 
\begin{align}
\begin{split}
\al_d^- &\defeq \frac{2}{d} \big( \log d - \log\log d \big) ;\\
\al_d^+ &\defeq \frac{2}{d} \log d
\label{eq:al+-}
\end{split}
\end{align}
\noindent (a) If $d \geq 66$, then $f_{d,\al}$ can have at most two local maxima.\\
(b) Moreover, if $d \geq 909$, then $f_{d,\al}$ has a (second) local maximum at some 
\[ \be_2^{\max}>\al- \frac{(\log d)^2}{d^{3/2}} .\]
\end{lemma}
\begin{proof}
Throughout this proof we think of $d,\al$ as fixed numbers satisfying the conditions of the lemma. Let
\begin{equation} \label{eq:be+-}
\be_- \defeq \frac{2}{d}  
\quad \text{and} \quad 
\be_+ \defeq \al - \frac{3}{d} 
\end{equation}
We claim that 
\begin{enumerate}[(i)]
\item $f''_{d,\al}(\be)>0$ on $[\be_-,\be_+]$;
\item the third derivative $f'''_{d,\al}$ is positive on $(0,\be_-)$ and negative on $(\be_+,\al)$.
\end{enumerate}

Part (a) follows quickly from these claims. Indeed, due to (i), $f''_{d,\al}$ has no root on $[\be_-,\be_+]$, while (ii) shows that it can have at most one root on each of the intervals $(0,\be_-]$ and $[\be_+,\al)$. Therefore, $f''_{d,\al}$ has at most two roots on $(0,\al)$, implying that $f'_{d,\al}$
can have at most three roots on $(0,\al)$. That is, $f_{d,\al}$ has at most three stationary points and hence at most two local maxima, as claimed.

We begin by noting that the first term in \eqref{eq:f_derder} is positive, while the second term is always greater than $(d-1)/(1-\al)$, hence we get the following lower bound: 
\[ f''_{d,\al}(\be) 
\geq g_{d,\al}(\be) \defeq
\frac{d-1}{1-\al}-\frac{2}{\al-\be}-\frac{1}{\be} .\]
Therefore, to show (i), it suffices to prove that $g_{d,\al}(\be)>0$ on the interval $[\be_-,\be_+]$. Since $g_{d,\al}(\be)$ is concave in $\be$, it is enough to prove positivity at the endpoints of this interval. Substituting $\be=\be_-=2/d$ we get 
\[ g_{d,\al}(\be_-) 
= \frac{d-1}{1-\al}-\frac{2}{\al-\frac{2}{d}}- \frac{d}{2} .\]
The expression on the right-hand side is clearly monotone increasing in $\al$, and hence we may consider the worst case scenario $\al=\al_d^-$. It is easy to see that, asymptotically, $g_{d,\al_d^-}(\be_-)$ is $(1/2-o_d(1)) d$ and it can be checked to be positive for any $d \geq 65$.

As for $\be=\be_+=\al-3/d$, we have
\[ g_{d,\al}(\be_+) 
= \frac{d-1}{1-\al}-\frac{2}{3}d- \frac{1}{\al-\frac{3}{d}} ,\]
which is asymptotically $(1/3-o_d(1)) d$. Again, we may consider the worst case $\al=\al_d^-$, when the right-hand side can be checked to be positive for any $d \geq 63$.

For the second claim (ii), recall formula \eqref{eq:f_derderder} expressing the third derivative $f'''_{d,\al}$. It is easy to see that the function $G_\al$ is monotone decreasing in $\be$, and hence $G'_\al$ is always negative. It also follows that 
\[ G_\al(\be) \geq G_\al(\al) 
= 2 \bigg( \frac12 - \al \bigg)^2 
\geq 2 \bigg( \frac12 - \al_d^+ \bigg)^2 
> \frac 14,\]
where the last inequality holds for all $d \geq 55$. We conclude that for $\be \leq \be_- = 2/d$ we have 
\begin{align*}
f'''_{d,\al}(\be) 
&> -4d + 0 
- \frac{d-1}{(1-2\al+\be)^2}
- \frac{2}{(\al-\be)^2}
+ \frac 1 {\be^2} \\
&\geq -4d 
- \frac{d-1}{(1-2\al_d^{+})^2}
- \frac{2}{\big( \al_d^{-} - \frac{2}{d} \big)^2}
+ \frac{d^2}{4} , 
\end{align*}
which is positive for all $d \geq 66$. It follows that $f'''_{d,\al}(\be)>0$ for all $0<\be<\be_-=2/d$.

Furthermore, it can be seen easily that the total contribution of the first two terms of \eqref{eq:f_derderder} is negative. Since the third term is also negative, we obtain 
\[ f'''_{d,\al}(\be) 
< -\frac{2}{(\al-\be)^2}+\frac 1 {\be^2} < 0 
\text{ provided that } \be > \big( \sqrt{2} -1 \big) \al .\]
We conclude that $f'''_{d,\al}(\be)<0$ for every $\be \geq \be_+=\al-3/d$ provided that $3/d < \big( 2-\sqrt{2} \big) \al$, which holds for every $\al \geq \al_d^-$ provided that $d \geq 51$. 

To prove (b), we first need to lower-bound the first derivative. We claim that for $\be \geq \al^2$ we have 
\[ 1/2-\al +\sqrt{(1/2-\al)^2+(\al-\be)^2} 
\leq 1-2\al+\al^2 .\]
Indeed, the left-hand side is clearly monotone decreasing in $\be$, and at $\be=\al^2$ it is equal to the right-hand side due to the identity \eqref{eq:identity}. Combining this with \eqref{eq:f_der} we get 
\[ f'_{d,\al}(\be) 
\geq -d \log(1-2\al+\al^2)+d \log(1-2\al+\be)
+2\log(\al-\be) 
\underbrace{-\log \be}_{\geq -\log \al} ,\]
where the sum of the first two terms can be further bounded, using \eqref{eq:log_upp_bnd}, as follows:
\begin{align*}
-d \log\left( \frac{1-2\al+\al^2}{1-2\al+\be} \right)
&= -d \log\left( 1 - \frac{\be-\al^2}{1-2\al+\be} \right) 
\geq d\,\frac{\be-\al^2}{1-2\al+\be} \\
&\geq d\,\frac{\al-\al^2+(\be-\al)}{1-\al} 
= d\al - d\,\frac{\al-\be}{1-\al} .
\end{align*} 
In conclusion:  
\begin{equation*} 
f'_{d,\al}(\be) 
\geq d\al - d\,\frac{\al-\be}{1-\al} + 2 \log(\al-\be)-\log \al
\quad \text{for any } \be \in [\al^2,\al) .
\end{equation*}
Setting $\be = \al-\de$ for some $\de<\al/2$ we get 
\[ f'_{d,\al}(\be) 
\geq d\al - \frac{d\de}{1-\al} + 2\log \de - \log \al .\]
In particular, for $\de^\star \defeq \frac{(\log d)^2}{d^{3/2}}$ this gives the following bound at $\be^\star \defeq \al-\de^\star$: 
\[ f'_{d,\al}(\be^\star) 
> d\al  - \frac{(\log d)^2/\sqrt{d}}{1-\al} 
- \log \al - 3\log d + 4 \log\log d .\]
This expression can be easily seen to be monotone increasing in $\al$. Therefore, we get the smallest value if we set $\al= \al_d^-$:
\[ 2(\log \log d) - \log\big( \log d 
- \log\log d \big) - \log 2 
- \frac{(\log d)^2}{\sqrt{d}(1-\al_d^-)} .\]
This turns out to be positive whenever $d \geq 909$. It follows that $f'_{d,\al}(\be^\star)>0$. Since $\lim_{\be\to\al} f'_{d,\al}(\be)=-\infty$, there must be a local maximum $\be_2^{\max}$ somewhere in $(\be^\star,\al)$. The proof is complete.
\end{proof}
%

\section{Bounding the second local maximum} \label{sec:second_max}

\newcommand\Coefd{\frac{2}{3}}
\newcommand\Coef{\frac{1}{3}}

Our goal in this section to prove an upper bound on $f_{d,\al}(\be_2^{\max})$, that is, on the second local maximum of $f_{d,\al}$, and then prove Theorem~\ref{thm:vs_alfm}. By Lemma~\ref{lem:local_maxima}(b) we already know that $\be_2^{\max}$ is near $\al$. So we start with estimating $f_{d,\al}(\be)$ around $\be=\al$. In fact, we will bound the difference $f_{d,\al}(\be) -f_{d,\al}(\al)$. The reason for this is that this difference is not so sensitive to the specific value of $\al$: it ``mainly depends'' on $\de = \al-\be$.
\begin{claim} \label{cl:f_diff}
Let $d \geq 318$ and $\al \leq (2\log d)/d$. Suppose that $\be=\al-\de$ for some $0\leq \de \leq \frac{3}{2d}$. Then 
\begin{equation*}
f_{d,\al}(\be) - \underbrace{f_{d,\al}(\al)}_{=\varphi_d(\al)} 
\leq 2 h(\de) 
+ \big( -d\al + \log(\al) + 2 \big) \de 
+ \Coefd d \de^2 .
\end{equation*}
\end{claim}
The proof, which is fairly technical, will be given in Section~\ref{sec:technical}. Here we deduce a corollary which bounds the maximum of $f_{d,\al}$ over the interval $[\al-\frac{3}{2d},\al]$. 
\begin{corollary} \label{cor:max_diff}
Let $d \geq 318$ and $(\log d + \log\log d + 2)/d \leq \al \leq (2\log d)/d$. Setting
\begin{equation}
c \defeq \frac{-d\al + \log(\al)}{2}+1 
\quad \text{and} \quad
C \defeq \Coef d,
\label{eq:cC}
\end{equation}
we define the function $\De_{d}(\al)$ as follows:
\begin{align}
\begin{split}
\De_{d}(\al) 
& \defeq 2\exp(c-1) \bigg( 1+ C \exp(c+1) \bigg) \\
&= 2\exp(c-1)+2C\exp(2c) \\
&= 2 \exp(-d\al/2) \sqrt{\al} 
+ \Coefd \exp(-d\al+2) d\al . 
\label{eq:Delta}
\end{split}    
\end{align}
Then 
\begin{equation*}
\max_{\be \in [\al-\frac{3}{2d},\al]} f_{d,\al}(\be) 
\leq \varphi_d(\al) + \De_{d}(\al) .
\end{equation*}
\end{corollary}
\begin{proof}
Due to the assumed (lower and upper) bounds on $\al$, we have 
\[ -d\al \leq -\log d - \log\log d - 2 
\quad \text{and} \quad 
\log \al \leq -\log d + \log\log d + \log 2 .\]
It follows that 
\[ c \leq -\log d + \frac{\log 2}{2} 
\text{, and hence } 
\exp(-c) \geq \frac{d}{\sqrt{2}} > \Coefd d = 2C .\]
Consequently, the condition of Claim~\ref{cl:hx+cx+Cxx} is satisfied so we can apply it to the function $g(\de) \defeq h(\de)+c\de + C \de^2$.
Therefore, 
\[ \max_{0 \leq \de \leq \frac{3}{2d}} g(\de) 
\leq \exp(c-1)+C\exp(2c) 
= \frac{\De_d(\al)}{2} .\]
Therefore, using Claim~\ref{cl:f_diff}, we conclude that  
\[ \max_{0 \leq \de \leq \frac{3}{2d}} 
f_{d,\al}(\al-\de) - \varphi_d(\al) \leq 
\max_{0 \leq \de \leq \frac{3}{2d}} 2g(\de) 
\leq \De_d(\al) ,\]
and the proof is complete.
\end{proof}
Using this corollary along with the fact that the second local maximum $\be_2^{\max}$ lies in the interval $[\al-\frac{3}{2d},\al]$ for large enough $d$, we get the following sufficient conditions for the second moment method.
\begin{lemma} \label{lem:sm_conds}
Suppose that 
$d \geq 4400$ and $\al \geq \al_d^- \defeq \frac{2}{d} (\log d-\log\log d)$.

\medskip

\noindent (a) \, If $\De_{d}(\al) \leq \varphi_d(\al)$, then the second moment condition \eqref{eq:the_cond} of Theorem~\ref{thm:Markovian_indep} is satisfied (i.e., $f_{d,\al}$ attains its maximum at $\al^2$), and hence $\alsm_d \geq \al$. 

\medskip

\noindent (b) \,  Suppose that for a positive $\eps$ we have 
\[ \al \leq \alfm_d - \eps 
\quad \text{and} \quad 
\De_{d}(\al) \leq -\varphi'_d(\al) \eps. \]
Then $\alsm_d \geq \alfm_d - \eps$.
\end{lemma}
\begin{proof}
The condition in (a) implies that $\varphi_d(\al)>0$, and hence $\al < \alfm_d < (2 \log d)/d$. Also, one can check that for $d \geq 4400$ we have 
%
\[ \al_d^- > \frac{1}{d}(\log d + \log\log d + 2) ,\]
and hence Corollary~\ref{cor:max_diff} can be applied. It follows that $f_{d,\al}$ is at most 
\[ \De_{d}(\al) + \varphi_d(\al)
\leq 2 \varphi_d(\al) 
= f_{d,\al}(\al^2) \]
over the interval $[\al-\frac{3}{2d},\al]$. According to Claim~\ref{cl:al_sq}, there is a local maximum at $\be_1^{\max}=\al^2$, while the two statements of Lemma~\ref{lem:local_maxima} yield that $f_{d,\al}$ has exactly one other local maximum $\be_2^{\max}$, which is located inside the interval $[\al-{(\log d)^2}/{d^{3/2}},\al]$. Since 
\[ \frac{(\log d)^2}{d^{3/2}} < \frac{3}{2d} 
\quad (d \geq 1029) ,\]
we get $\be_2^{\max} \in [\al-\frac{3}{2d},\al]$. We conclude that $f_{d,\al}(\be_2^{\max}) \leq f_{d,\al}(\al^2)$. Therefore, $f_{d,\al}$ attains its global maximum at $\al^2$, as claimed. 

To see (b) let $\alh \defeq \alfm_d-\eps \geq \al$. Using $\varphi_d(\alfm_d)=0$, we get 
\[ -\varphi_d(\alh) 
= \varphi_d(\alfm_d) -\varphi_d(\alh) 
= \int_{\alh}^{\alfm_d} \varphi'_d .\]
Since $\varphi'_d$ is monotone decreasing (see Claim~\ref{cl:phi_der} in the Appendix),  we can upper bound this as follows:
\[ -\varphi_d(\alh) 
\leq \underbrace{(\alfm_d-\alh)}_{=\eps} \,
\underbrace{\varphi'_d(\alh)}_{\leq \varphi'_d(\al)} 
\leq \eps \varphi'_d(\al) .\]
We will see in a moment that $\De_d$ is monotone decreasing. It follows that 
\[ \varphi_d(\alh) 
\geq -\eps \varphi'_d(\al) 
\geq \De_d(\al) \geq \De_d(\alh) ,\]
so we can use (a) for $\alh$ to get $\alsm_d \geq \alh = \alfm_d - \eps$.

To see that $\De_d(\al)=2\exp(c-1)+2C\exp(2c)$ is monotone decreasing in $\al$, recall \eqref{eq:cC} and notice that $C$ is a constant (for a fixed $d$) and $c$ is monotone decreasing for $\al > 1/d$. 
\end{proof}

The significance of Lemma~\ref{lem:sm_conds}(a) is that it provides an explicit condition involving reasonably simple, concrete functions ($\varphi_d$ and $\De_{d}$) so one can easily and very quickly check the condition with a computer for a specific pair $d, \al$. As for part (b), we will use it next to prove Theorem~\ref{thm:vs_alfm} for large degrees $d \geq d_0$.
\begin{proof}[Proof of Theorem~\ref{thm:vs_alfm}]
We need to consider three regimes for the degree $d$.

In the intermediate regime $2\cdot 10^5 \leq d \leq 2\cdot 10^6$ we used a computer to verify the condition $\De_{d}(\al) \leq \varphi_d(\al)$ of Lemma~\ref{lem:sm_conds}(a) for $\al = \alfm_d-\eps_d$ (which is indeed larger than $\al_d^-$, as required by the lemma). This verification can be done very quickly and with high precision by substituting the particular values of $d$ and $\al$ into the formulas \eqref{eq:Delta} and \eqref{eq:phi}, expressing $\De_{d}(\al)$ and $\varphi_d(\al)$, respectively. For instance, for $d=2\cdot 10^5$ we get 
\begin{align*}
\De_{d}(\al) & \approx   9.071810568001 \cdot 10^{-7};\\
\varphi_d(\al) & \approx 9.091179026828 \cdot 10^{-7}.
\end{align*}
Since $\De_{d}(\al) \leq \varphi_d(\al)$, it follows from the lemma that $\alsm_d \geq \al$. To get an idea how sharp the obtained inequality is, we computed the actual numerical values (for $d=2\cdot 10^5$): 
\begin{align*}
\al = \alfm_d-\eps_d &\approx 0.00010182262 ;\\
\alsm_d &\approx              0.00010183455;\\
\alfm_d &\approx              0.00010190391.
\end{align*}

In the ``small'' regime $d<2\cdot 10^5$ the condition $\De_{d}(\al) \leq \varphi_d(\al)$ is not satisfied.
We can still use a computer in this regime where we observe the same picture and we can show that the first local maximum at $\be_1^{\max}=\al^2$ is larger than the second at $\be_2^{\max}$: $f_{d,\al}(\be_1^{\max})>f_{d,\al}(\be_2^{\max})$. 
More precisely, we can get a rigorous numerical upper bound on $f_{d,\al}(\be_2^{\max})$ as follows. Given a concave function $g(\xi)$, one can use a binary search to find a small interval $[\xi_1,\xi_2]$ such that the function value at the midpoint $\xi_{\textrm{mid}}=(\xi_1+\xi_2)/2$ is larger than both $g(\xi_1)$ and $g(\xi_2)$. Then the maximum is certainly attained inside this interval and we have 
\[ \max g 
\leq \max_{[\xi_1,\xi_2]} g 
\leq 2g(\xi_{\textrm{mid}})-\min \big( g(\xi_1), g(\xi_2) \big) .\] 
Since $f_{d,\al}$ is concave in some interval $(\be_2^{\max}-\eps,\al)$, we can use the method sketched above to get a rigorous numerical bound on $f_{d,\al}(\be_2^{\max})$ and verify that the second moment condition \eqref{eq:the_cond} holds with $\al=\alfm_d - \eps_d$ for every $d \leq 2 \cdot 10^5$.

So it remains to prove the theorem for large degrees $d>d_0=2 \cdot 10^6$. First, for the sake of brevity, we introduce the following parametrization of the density $\al$: for a parameter $a$ let
\begin{equation} \label{eq:al_a}
\al_{d,a} \defeq 
\frac{2}{d}\bigg( \log d - \log\log d 
+ \ka + \frac{(\log\log d) - a}{\log d} \bigg) 
\text{, where } \ka \defeq 1-\log 2 \text{ is a constant.}
\end{equation} 
Note that, for fixed $d$, $\al_{d,a}$ is monotone decreasing in $a$. Furthermore, if we set $a=\ka$, then the density $\al_{d,\ka}$ is below the first moment bound $\alfm_d$ (for sufficiently large $d$). More precisely, Theorem~\ref{thm:alfm_bnd}(c) in the Appendix shows that 
\begin{equation} \label{eq:ineq2}
\al_{d,\ka} < \alfm_d 
\quad (d \geq 2 \cdot 10^6) .
\end{equation}%
Our proof strategy is the following. We fix a constant $a>\ka=1-\log 2 \approx 0.3$ so that $\al_{d,a} < \al_{d,\ka} < \alfm_d$. The argument would work for any $a>\ka$ but in order to keep $d_0$ reasonably small, we need that $a$ is neither too large nor too close to $\ka$. We will use the following specific choice: $a=\ka+0.2 \approx 0.5$.

We aim to prove that the following inequalities hold simultaneously for $d \geq d_0$: 
\begin{align}
\De_d(\al_{d,a}) &< -\varphi'_d(\al_{d,a}) \eps_d ;
\label{eq:ineq1} \\
\eps_d &< \al_{d,\ka}-\al_{d,a} 
= \frac{2(a-\ka)}{d(\log d)} .
\label{eq:ineq3} 
\end{align}
The fairly technical proof of \eqref{eq:ineq1} will be given in Section~\ref{sec:technical}, see Claim~\ref{cl:De_per_phider}. As for \eqref{eq:ineq3}, we have $a-\ka=0.2$, so it suffices to have:
\[ \underbrace{\frac{4\sqrt{2}}{e} \frac{(\log d)^{1/2}}{d^{3/2}}}_{\eps_d} < \frac{0.4}{d (\log d)} 
\Longleftrightarrow 
\frac{10\sqrt{2}}{e} (\log d)^{3/2} < \sqrt{d} 
\Longleftrightarrow 
\frac{200}{e^2} (\log d)^{3} < d ,\]
which is true for all $d \geq 3 \cdot 10^4$.

Finally we show why (\ref{eq:ineq1}--\ref{eq:ineq3}) imply $\alsm_d \geq \alfm_d - \eps_d$. Combining \eqref{eq:ineq2} and \eqref{eq:ineq3} shows that $\al_{d,a}< \al_{d,\ka}-\eps_d <\alfm_d-\eps_d$. This, along with \eqref{eq:ineq1}, means that we can apply Lemma~\ref{lem:sm_conds}(b) with $\al=\al_{d,a}$ and $\eps=\eps_d$ to conclude that $\alsm_d \geq \alfm_d-\eps_d$, as claimed.
\end{proof}
%

\section{Technical proofs} \label{sec:technical}

\begin{proof}[Proof of Claim~\ref{cl:f_diff}]
For the sake of neater formulas, we will simply \underline{write $\ga$ for $\Ga(\be)$} throughout the proof, where $\be=\al-\de$ as in the statement. With this mind, we can write $f_{d,\al}(\be)$ as follows:
\begin{align*}
f_{d,\al}(\be) 
&= \frac{d}{2} \bigg( 2h(\be)+2h(\ga)
+4 h(\al-\be-\ga) + h(1-4\al+2\be+2\ga) \bigg) \\
&- (d-1) \bigg( h(\be)+ 2h(\al-\be) 
+ h(1-2\al+\be) \bigg) \\
&= h(\al-\de) + d h(\ga) 
+ 2d \bigg( h(\de-\ga) - h(\de) \bigg) + 2h(\de) \\
&+ \frac{d}{2} h\big(1-2\al-2(\de-\ga)\big) 
- (d-1) h(1-\al-\de) ,
\end{align*}  
while \eqref{eq:Ga_second} translates to 
\begin{equation} \label{eq:ga_bounds}
\ga= \frac{\de^2}{\sqrt{\big(\frac12-\al\big)^2+\de^2}+\big(\frac12-\al\big)} .
\end{equation}
Note that we have the following bounds for $\ga$:
\begin{equation} \label{eq:ga_bnd}
\de^2 
\leq \frac{1}{1-\al} \de^2
\leq \ga 
\leq \frac{1}{1-2\al} \de^2 .
\end{equation}

Now we turn our attention to the difference in question: fix $d,\al$ and let 
\begin{equation} \label{eq:B_de}
B_\de \defeq f_{d,\al}(\be)-f_{d,\al}(\al) . 
\end{equation}
For brevity we use the notation $\de' \defeq \de-\ga$. Then 
\begin{align*}
B_\de 
&= \bigg( h(\al-\de) - h(\al) \bigg) + d h(\ga) 
+ 2d \bigg( h(\de') - h(\de) \bigg) + 2h(\de) \\
&+ \frac{d}{2} \bigg( h(1-2\al-2\de') - h(1-2\al) \bigg)
- (d-1) \bigg( h(1-\al-\de) - h(1-\al) \bigg).
\end{align*}  
Using \eqref{eq:h_diff1} with $x=\al$ and $\eps=\de$ we get 
\begin{equation} \label{eq:bnd1}
h(\al-\de) - h(\al) \leq -\de h'(\al) 
= \de \big( 1+\log \al \big) .
\end{equation}
Using \eqref{eq:h_diff1} again, this time with $x=\de$ and $\eps=\ga$ so that $x-\eps=\de'$, we get 
\[ h(\de') - h(\de) \leq -\ga h'(\de) 
= \ga \big( 1+ \log \de \big) ,\]
and hence 
\begin{equation} \label{eq:bnd2}
2d\bigg( h(\de') - h(\de) \bigg) + d h(\ga) 
\leq d \ga \big(2 + 2 \log \de  - \log \ga \big) 
\leq 2d \ga ,
\end{equation}
where we used $2 \log \de \leq \log \ga$ as $\de^2 \leq \ga$ according to \eqref{eq:ga_bnd}.

Recall that we work under the assumptions $\al\leq 2(\log d)/d$ and $\de \leq \frac{3}{2d}$. Then, according to \eqref{eq:ga_bnd}, we have 
\[ \frac{\ga}{\de} 
\leq \frac{\de}{1-2\al}
\leq \frac{3/(2d)}{1-4(\log d)/d} 
= \frac{3/2}{d-4(\log d)} < 0.37 
\text{ for any } d \geq 15. \]
It follows that 
\begin{equation} \label{eq:deperde}
\frac{\de'}{\de} = 1 - \frac{\ga}{\de} 
> 1- 0.37 
> \frac{1}{\sqrt[3]{4}} .
\end{equation}
Using \eqref{eq:h_diff2_simple} with $x=\al$ and $\eps=\de$, then substituting $x=2 \al$ and $\eps=2 \de'$ into the same formula, we get the following: 
\begin{align*} 
H_1 &\defeq h(1-\al-\de) - h(1-\al) 
= \de - \al\de - \frac{\de^2}{2}  
- \underbrace{\sum_{k=3}^\infty 
\frac{(\al+\de)^k-\al^k}{k(k-1)} }_{S_1 \defeq} ;\\
H_2 &\defeq h(1-2\al-2\de') - h(1-2\al) 
= 2 \de' - 4 \al \de' - 2 (\de')^2 
- \underbrace{\sum_{k=3}^\infty 
2^k \frac{(\al+\de')^k-\al^k}{k(k-1)}}_{S_2 \defeq} .
\end{align*}
Therefore, using $\de'=\de-\ga$, 
\begin{align} 
\begin{split}
\frac{d}{2} H_2 - (d-1)H_1
&= d\de'-2d\al\de' - d (\de')^2 - \frac{d}{2} S_2\\
&- (d-1) \de 
+ \underbrace{(d-1)\al \de + (d-1)\de^2/2 + (d-1) S_1}_
{\leq d\al \de + d\de^2/2 + d S_1} \\
&\leq \big( -d\al + 1 \big) \de + 
d\left( -\ga -(\de')^2 + \frac{\de^2}{2} + 2\al\ga \right) 
+ d\underbrace{ \left( S_1 - \frac{S_2}{2} \right)}_{\leq 0} ,
\label{eq:bnd3}
\end{split}
\end{align}
where the last term is negative due to the following observation: after expansion, we have 
\[ S_1 = \sum_{k=3}^\infty \frac{(\al+\de)^k-\al^k}{k(k-1)} 
= \sum_{k=3}^\infty \sum_{i=1}^k a_{k,i} \al^{k-i} \de^i 
\text{ and } 
\frac{S_2}{2} = \sum_{k=3}^\infty \sum_{i=1}^k 
2^{k-1} a_{k,i} \al^{k-i} (\de')^i
\]
for some positive coefficients $a_{k,i}$. Due to \eqref{eq:deperde}, $\de^i \leq 2^{k-1} (\de')^i$ holds for any $i \leq k$ and $k \geq 3$, and $S_1 \leq S_2/2$ follows.  

Summing up the bounds \eqref{eq:bnd1}, \eqref{eq:bnd2}, \eqref{eq:bnd3} and adding $2 h(\de)$:
\[ B_\de 
\leq 2 h(\de) 
+ \big( -d\al + \log(\al) + 2 \big) \de 
+ d\bigg( 
\underbrace{\ga -(\de')^2 + \frac{\de^2}{2} + 2\al\ga}_
{R \defeq} \bigg) .\]
To complete the proof, it suffices to show that $R \leq \Coefd \de^2$ for the last term on the right-hand side. First note that 
\[ (\de')^2 = (\de-\ga)^2 \geq \de^2 - 2\de \ga .\]
Therefore, using \eqref{eq:ga_bnd}, 
\[ R \leq (1+2\al+2\de) \ga - \de^2/2 
\leq \left( \frac{1+2\al+2\de}{1-2\al} 
- \frac{1}{2} \right) \de^2 .\]
It remains to show that the coefficient of $\de^2$ is less than $\Coefd$. Note that this coefficient is monotone increasing both in $\al$ and in $\de$, so we may set them to their maximal values, that is, $\al=(2 \log d)/d$ and $\de=3/(2d)$, in which case the coefficient becomes  
\[ \frac{1+2\al+2\de}{1-2\al} - \frac{1}{2} 
= \frac{d/2+6 \log d + 3}{d-4 \log d} ,\]
which is indeed less than $\Coefd$ provided that $d \geq 318$.
\end{proof}

Finally, we prove the last missing ingredient: inequality \eqref{eq:ineq1}.
\begin{claim} \label{cl:De_per_phider}
Let $a=\ka+0.2$ and $d \geq d_0 = 2 \cdot 10^6$. Then 
\[ \De_d(\al_{d,a}) < -\varphi'_d(\al_{d,a}) \eps_d .\]
\end{claim}
\begin{proof}
Under the assumption $d \geq 2 \cdot 10^6$ we also have
\begin{align}
\begin{split}
\log d &> 14.5 ;\\
\log \log d &> 2.67 ;\\
\frac{\log \log d}{\log d} &< 0.185 < \frac 15 . 
\label{eq:d_bnds}
\end{split}
\end{align}
For brevity we will write $\al$ for $\al_{d,a}$. Recall that we defined $\al_{d,a}$ in \eqref{eq:al_a}, which is a monotone decreasing expression in $a$. Since $a> \ka$, we have 
\begin{align*}
\al &= \al_{d,a} < \al_{d,\ka} 
= \frac{2 \log d}{d} 
\left( 1 - \frac{(\log\log d)-\ka}{\log d} 
+ \frac{(\log\log d) - \ka}{(\log d)^2} \right) \\
&= \frac{2 \log d}{d} 
\left( 1 - x + \frac{x^2}{(\log\log d) - \ka} \right) 
\text{, where } x \defeq \frac{(\log\log d)-\ka}{\log d} .
\end{align*} 
Since $(\log \log d)-\ka> 2.36$ by \eqref{eq:d_bnds}, we get that  
\[ 1-x+\frac{x^2}{(\log \log d)-\ka} 
< 1 - x + \frac{x^2}{2.36} < \exp(-x) ,\]
where the second inequality can be checked to hold for any $x<1/2$. This is the case now because $x<(\log\log d)/(\log d) <0.185$ by \eqref{eq:d_bnds}. We conclude that 
\begin{equation}
\label{eq:Tbound1}
    \al < \frac{2 \log d}{d} 
\exp \underbrace{ \left(- \frac{\log\log d}{\log d} 
+ \frac{\ka}{\log d} \right)}_{T_1 \defeq}  
\end{equation}
Next we set $c=-d\al/2+(\log \al)/2+1$ and $C = \Coef d$ as in \eqref{eq:cC} of Corollary~\ref{cor:max_diff}. In the current setting we have $\al=\al_{d,a}$ so  
\begin{align}
\begin{split}
c-1 &= -\log d + (\log\log d) - \ka 
- \frac{(\log\log d) - a}{\log d} + \frac{\log \al}{2} 
\text{, and hence} \\
\exp(c-1) &= \frac{2 \log d}{e d} \sqrt{\al} 
\,\exp\underbrace{ \left( 
\frac{a-\log\log d}{\log d} \right) }_{T_2 \defeq} \\
<& \frac{2\sqrt{2}}{e} \left( \frac{\log d}{d} \right)^{3/2} 
\exp\left( \frac{T_1}{2}+T_2 \right) \quad 
\text{ by \eqref{eq:Tbound1}.}
\label{eq:Tbound2}
\end{split}
\end{align} 
Using $\al < (2 \log d)/d$ and $1+x \leq \exp(x)$, from \eqref{eq:Delta} we get 
\begin{align}
\begin{split}\label{eq:Tbound3}
\De_d(\al) &= 2\exp(c-1) \bigg( 1+ C \exp(c+1) \bigg) 
\leq 2\exp(c-1) \bigg( 1+ \Coefd e (\log d) \sqrt{\al} \bigg) \\
&\leq 2\exp(c-1) \bigg( 1+ \Coefd \sqrt{2} e 
\frac{(\log d)^{3/2}}{d^{1/2}} \bigg) 
\leq 2\exp(c-1) \exp \underbrace{\bigg( \Coefd \sqrt{2} e 
\frac{(\log d)^{3/2}}{d^{1/2}} \bigg)}_{T_3 \defeq} .
\end{split}
\end{align}
Furthermore, one can show easily that $\varphi'_d(\al) < -\log \al - d\al$ (see Claim~\ref{cl:phi_der} in the Appendix). Since $a<2.67<\log \log d$, we have $\al > \al_{d,\log\log d}$. One can bound $\log(\al)$ easily (see \eqref{eq:log_al_minus} in Claim~\ref{cl:log_al_minus} of the Appendix) and get that   
\begin{align*}
-\varphi'_d(\al) &> \log \al + d\al \\
&> -\log d + \log\log d + \log 2 
-\frac{\log\log d}{\log d} -\frac{1/4}{\log d}\\
&\quad + 2\bigg( \log d - \log\log d + \ka 
+ \frac{(\log\log d) - a}{\log d} \bigg) \\
&= \log d - \log\log d + 2 - \log(2)
+ \underbrace{\frac{(\log\log d)-2a-1/4}{\log d}}_{>0} , 
\end{align*}
where the last term is positive because $\log\log d > 2.67 \geq 2a+1/4$. Now let 
\[ x \defeq \frac{(\log\log d)-2+\log 2}{\log d} 
< \frac{\log\log d}{\log d} < 0.185 .\]
By \eqref{eq:log_low_bnd1} we have $1-x > \exp(-x-x^2)$. Therefore 
\begin{align}
\begin{split}
\label{eq:Tbound4}
-\varphi'_d(\al) &> (\log d) (1-x) 
> (\log d)\exp(-x-x^2) \\
&> (\log d) \exp \underbrace{\left( 
\frac{-(\log\log d)+2-\log 2}{\log d} 
- \frac{(\log\log d)^2}{(\log d)^2} \right) }_{T_4 \defeq} .
\end{split}
\end{align}
Combining (\ref{eq:Tbound1}--\ref{eq:Tbound4}) yields the following bound: 
\begin{align*}
\frac{\De_d(\al)}{-\varphi'_d(\al)} &< \eps_d 
\cdot \exp(T) , \quad \text{where } \\
\eps_d &= \frac{4\sqrt{2}}{e} 
\frac{(\log d)^{1/2}}{d^{3/2}} ;\\
T & \defeq \frac{T_1}{2}+T_2+T_3-T_4 .
\end{align*}
It suffices to show that $T<0$. To this end, we write $T$ out explicitly:
\begin{align*}
T &= \frac12 \left( - \frac{\log\log d}{\log d} + \frac{\ka}{\log d} \right) 
+ \frac{a-\log\log d}{\log d} 
+ \Coefd \sqrt{2} e 
\frac{(\log d)^{3/2}}{d^{1/2}} \\
&+ \frac{(\log\log d)-2+\log 2}{\log d} 
+\frac{(\log\log d)^2}{(\log d)^2} \\
&= \frac{-\frac12 (\log\log d) + a + \log(2) + \ka/2 - 2 
+ (\log\log d)^2/(\log d)}{\log d}
+ \Coefd \sqrt{2} e 
\frac{(\log d)^{3/2}}{d^{1/2}} .
\end{align*}
It is easy to see that, for any fixed $a$, $T$ is eventually negative (i.e., for sufficiently large $d$). For the specific choice $a=\ka+0.2$, $T$ can be checked to be negative if $d \geq 1.9\cdot 10^6$, and the proof is complete. 
\end{proof}

\bibliographystyle{plain}
\bibliography{refs}

\newpage
\section{Appendix: the first moment bound} \label{sec:fm_bound}

The following function (and its unique positive root) was first introduced by Bollob\'as in \cite{bollobas1981independence} for studying the independence ratio of random regular graphs: for $\al \in [0,1/2]$ let 
\begin{align*}
\varphi_d(\al)  
&\defeq h(\al)+h(1-\al) + \frac{d}{2} 
\bigg( h(1-2\al) - 2h(1-\al) \bigg) \\
&= h(\al) + \frac{d}{2} h(1-2\al) - (d-1) h(1-\al) 
\text{, where } h(x) \defeq -x \log x .
\end{align*}
It expresses the ``first moment entropy'', i.e., the exponential rate of the expected number of independent sets of density $\al$ (for random $d$-regular graphs). In particular, if $\varphi_d(\al)<0$, then the probability that the graph contains any independent set of density $\al$ is exponentially small. It follows that the asymptotic independence ratio $\alstar_d$ is less than $\al$. In conclusion: 
\[ \alstar_d \leq \alfm_d 
\defeq \inf\big\{ \al \in [0,1/2] \, : \, \varphi_d(\al)<0 \big\} .\]
We call $\alfm_d$ the \emph{first moment (upper) bound}. We will see that $\varphi_d(\al)/\al$ is continuous and monotone decreasing, and hence has a unique root on $(0,1/2)$. It follows that $\alfm_d$ is the unique root of $\varphi_d$.

The asymptotic formula \eqref{eq:al_asymptotic} approximates $\alfm_d$ with a (non-explicit) error term $o(1/d)$. The Ding--Sly--Sun paper mentions (without a proof) the following finer asymptotic formula \cite[formula (4)]{ding2016maximum}:
\begin{equation*} 
\alfm_d = \frac{2}{d} \left( \log d - \log \log d + 1 - \log 2 
+ \Oc\left( \frac{\log \log d}{\log d} \right) \right) .
\end{equation*}
Our goal is to rigorously prove the following (even more accurate) estimates. 
\begin{theorem} \label{thm:alfm_bnd}
Let 
\begin{equation*} 
\ka \defeq \log(e/2) = 1-\log 2 \approx 0.3,    
\end{equation*}
and for a real parameter $a \geq 0$ define 
\begin{equation*} 
\al_{d,a} \defeq 
\frac{2}{d}\bigg( \log d - \log\log d 
+ \ka + \frac{(\log\log d) - a}{\log d} \bigg) .
\end{equation*} 
Then, for any fixed $a \geq 0$, there exists $d_0(a) \in \Nb$ such that 
\begin{align*}
\text{if } a &< \ka 
\text{, then } \alfm_d < \al_{d,a} \;
\text{ for } d \geq d_0(a) ; \tag{a}\\
\text{if } a &\geq \ka 
\text{, then } \alfm_d > \al_{d,a} \;
\text{ for } d \geq d_0(a) \tag{b}.
\end{align*}
In particular, for the case $a=\ka$ we will prove the following:
\[ \alfm_d > \al_{d,\ka} 
\; \text{ for } d \geq 2 \cdot 10^6 . \tag{c} \]
One can actually verify with a computer that (c) holds for any $d \geq 282118$.
\end{theorem}
\begin{remark}
Another way to phrase this theorem is the following. For a given $d$, we define $a_d$ to be the unique real number for which $\alfm_d=\al_{d,a_d}$. Then the following equivalence clearly holds:
\[ \alfm_d > \al_{d,a} \Longleftrightarrow a_d < a .\] 
So an equivalent formulation of (a) and (b) is that $a_d \to \ka$ in such a way that $a_d< \ka$ for large enough $d$, while (c) claims that $a_d< \ka$ ($d \geq 2 \cdot 10^6$).

In fact, the sequence $(a_d)_{d\geq 3}$ behaves as follows. First it grows from $a_3 \approx 0.77$ to its maximal value $a_{23} \approx 1.75$. Then it decreases, and it appears to take its minimal value $\approx 0.2489$ at around $d \approx 2^{64}$. After that, $a_d$ seems to be monotone increasing again, converging to $\ka \approx 0.30685$.
\end{remark}
\begin{remark}
It is straightforward to check with a computer if $\alfm_d$ is below or above $\al_{d,a}$ for specific $d$ and $a$: one simply needs to evaluate $\varphi_d(\al_{d,a})$ and check its sign (which can be done very quickly and with high precision).
\begin{itemize}
\item If $\varphi_d(\al_{d,a})>0$, then $\alfm_d>\al_{d,a}$.
\item If $\varphi_d(\al_{d,a})<0$, then $\alfm_d<\al_{d,a}$.
\end{itemize}
According to Theorem~\ref{thm:alfm_bnd}(c), for any $a \geq \ka$ we have 
\[ \alfm_d > \al_{d,\ka} \geq \al_{d,a} 
\quad (d \geq 2 \cdot 10^6) .\]
It follows that $d_0(a)$ can be chosen to be $2 \cdot 10^6$ in this case. Using a computer one can determine the precise threshold for a given $a \geq \ka$. For instance, we have 
\begin{align*}
d_0(\ka) &= 282118 ;\\
d_0(0.4) &= 15072 ;\\
d_0(1) &= 392 .
\end{align*}
As for the case $a<\ka$, computer checks suggest that $d_0(a)=3$ for any $0 \leq a \leq 0.24$, that is, $\alfm_d < \al_{d,0.24}$ holds for all degrees $d \geq 3$. When $a$ is just below $\ka$, however, $d_0(a)$ is gigantic. For instance, 
\[ d_0(0.28) \approx 2^{556} \approx 2.3 \cdot 10^{167} .\]
\end{remark}
\subsection{Derivative and related functions}
Using $h'(x)=-1-\log x$ we get the following for the derivative $\varphi'_d$:
\begin{align*}
\varphi'_d(\al) 
&= -\log(\al) + d \log(1-2\al) - (d-1) \log(1-\al) \\
&= -\log(\al) - (d+1) \al 
- \sum_{k \geq 2} \frac{(2^k-1)d+1}{k}\al^k ,
\end{align*}
where the second equality follows from the power series \eqref{eq:log_power_series} of $\log(1-x)$. The following is an immediate consequence of this expansion.
\begin{claim} \label{cl:phi_der}
$\varphi'_d$ is monotone decreasing on $(0,1/2)$ and satisfies 
\[ \varphi'_d(\al) < -\log \al - d\al . \]
Furthermore,
\[ \lim_{\al \to 0+ } \varphi'_d(\al) = +\infty 
\quad \text{and} \quad 
\lim_{\al \to \frac 12 - } \varphi'_d(\al) = -\infty .\]
\end{claim}

Often it will be more convenient to work with the following variant of $\varphi_d$ instead:
\begin{align}
\begin{split} \label{eq:psi}
\psi_d(\al) \defeq \frac{\varphi_d(\al)}{\al} 
&= -\log(\al) + \frac{h(1-\al)}{\al} 
+ \frac{d}{2} \frac{h(1-2\al) - 2h(1-\al)}{\al} \\
&= -\log(\al) + 1 - \frac{d}{2}\al 
- \sum_{k=1}^\infty \frac{\al^k}{k(k+1)} 
- d \sum_{k=2}^\infty \frac{2^k-1}{k(k+1)} \al^k ,
\end{split}
\end{align}
where we used \eqref{eq:h_ps} to get the expansion. This expansion clearly shows that $\psi_d(\al)$ is continuous and monotone decreasing on $(0,1/2)$. Moreover,
\[ \lim_{\al \to 0+} \psi_d(\al) = \infty 
\quad \text{and} \quad 
\psi_d(1/2) = -(d-2) \log 2 .\]
Consequently, $\psi_d$ (and hence $\varphi_d$) has a unique root on $(0,1/2)$. This root is the first moment bound $\alfm_d$.

Let us introduce the following notations for the main term and the remainder term in the expansion \eqref{eq:psi}:
\begin{align}
\psi_d(\al) &= \psimain_d(\al) - \psirem_d(\al) 
\text{, where} \nonumber\\
\psimain_d(\al) &\defeq -\log(\al) + 1 - \frac{d}{2}\al ;
\label{eq:psimain} \\
\psirem_d(\al) &\defeq 
\sum_{k=1}^\infty \frac{\al^k}{k(k+1)} 
+ d \sum_{k=2}^\infty \frac{2^k-1}{k(k+1)} \al^k .
\label{eq:psirem}
\end{align} 
Next we give bounds for the remainder term.
\begin{claim} \label{cl:remainder} 
We have the following bounds for the remainder term $\psimain_d(\al) - \psi_d(\al)$:
\[ \frac 12 \al + \frac{d}{2}\al^2 
< \; \psimain_d(\al) - \psi_d(\al) \;
< \frac 12 \al + \frac{d+1}{2}\al^2 
+ \frac 23 d \frac{\al^3}{1-2\al} .\]
\end{claim}
\begin{proof}
The lower bound follows immediately from \eqref{eq:psirem}.

As for the upper bound: on the one hand, we have 
\[ \sum_{k=1}^\infty \frac{\al^k}{k(k+1)} 
\leq \frac 12 \al
+ \sum_{k=2}^\infty \frac{\al^2}{k(k+1)} 
= \frac 12 \al + \frac 12 \al^2 ;\]
on the other hand, 
\[ \sum_{k=2}^\infty \frac{2^k-1}{k(k+1)} \al^k 
\leq \frac 12 \al^2 
+ \sum_{k=3}^\infty \frac{2^k}{12} \al^k 
= \frac 12 \al^2 
+ \frac{8 \al^3}{12} \sum_{k=0}^\infty (2\al)^k 
= \frac 12 \al^2 
+ \frac 23 \, \frac{\al^3}{1-2\al} .\]
These two inequalities clearly yield the claimed bound on the remainder term.
\end{proof}

The following is an immediate corollary.
\begin{lemma} \label{lem:implications}
We have the following implications:

\medskip

\noindent (a) \;
$\ds \psimain_d(\al) 
< \frac 12 \al + \frac{d}{2}\al^2 
\Longrightarrow 
\psi_d(\al) < 0 
\Longrightarrow
\alfm_d < \al$ ;

\medskip

\noindent (b) \;
$\ds \psimain_d(\al) 
> \frac 12 \al + \frac{d+1}{2}\al^2 
+ \frac 23 d \frac{\al^3}{1-2\al} 
\Longrightarrow 
\psi_d(\al) > 0 
\Longrightarrow
\alfm_d > \al$ .
\end{lemma}

One can quickly deduce the basic upper bound 
\begin{equation} \label{eq:al_basic_upper}
\alfm_d < \frac{2 \log d}{d} .
\end{equation}
Indeed, set $\al=\frac{2\log d}{d}$ and note that 
\[ \psimain_d(\al) = \log d - \log\log d - \log(2) + 1 - \log d = - \log\log d - \log(2) + 1 
< 0 \text{ for } d \geq 4 .\]
By Lemma~\ref{lem:implications}(a) we conclude that $\alfm_d < \al$.


\subsection{Logarithm of the density}
Our next goal is to give precise estimates for $\log \al_{d,a}$. First we rewrite \eqref{eq:al_a} as 
\begin{align}
\al_{d,a} &= \frac{2\log d}{d} (1-x+y) \text{, where } \notag \\
x &=\frac{(\log\log d)-\ka}{\log d} \text{ and} 
\label{eq:x} \\
y &=\frac{(\log\log d)-a}{(\log d)^2} .
\label{eq:y}
\end{align}
Then we have 
\begin{equation}
\log \al_{d,a} 
= -\log d + \log\log d + \log 2 + \log(1-x+y) .
\label{eq:log_al_a} 
\end{equation} 
The main term in $\log(1-x+y)$ is $-x$ so we introduce the following notation for the remainder term:
\begin{equation} \label{eq:gxy}
g(x,y) \defeq -x-\log(1-x+y)
= -y + \sum_{k=2}^\infty \frac{(x-y)^k}{k} ,
\end{equation}
where the expansion follows from \eqref{eq:log_power_series}. Then 
\begin{equation}
\log \al_{d,a} 
= -\log d + \log\log d + \log 2 
- \frac{(\log\log d)-\ka}{\log d} - g(x,y) .
\label{eq:log_al_a_g_version} 
\end{equation} 
Using \eqref{eq:al_a} and \eqref{eq:log_al_a_g_version} we can express $\psimain_d(\al_{d,a})$ as follows:
\begin{align}
\begin{split}
\psimain_d(\al_{d,a}) 
&=  -\log(\al_{d,a}) + 1 - \frac{d}{2} \al_{d,a} \\
&= \log d - \log\log d - \log(2) 
+ \frac{(\log\log d)-\ka}{\log d} + g(x,y) \\
&+ 1 - \left( \log d - \log\log d + \ka 
+ \frac{(\log\log d) - a}{\log d} \right)\\
&= \frac{a-\ka}{\log d} + g(x,y) .
\label{eq:Psi_g}
\end{split}
\end{align}
From \eqref{eq:x}, \eqref{eq:y}, and \eqref{eq:gxy} it is easy to see that 
\begin{equation} \label{eq:gxy_asymptotics}
g(x,y) = \left(\frac 12 - o_d(1) \right) 
\left( \frac{\log\log d}{\log d} \right)^2 .
\end{equation}
Note that this convergence is very slow and $g(x,y)$ is actually negative if $d$ is not large enough. Later we will prove the following lower bounds on $g(x,y)$.
\begin{claim} \label{cl:g_pos}
Suppose that $a \geq \ka$. Then:
\begin{align*}
\text{if } d &\geq 1.1 \cdot 10^6
\text{, then } g(x,y)>0 
\text{ and } \psimain_d(\al_{d,a}) > \frac{a-\ka}{\log d} ;\\
\text{if } d &\geq 2 \cdot 10^6
\text{, then } g(x,y)> \frac{x}{40 \log d} 
\text{ and } \psimain_d(\al_{d,a}) > \frac{a-\ka}{\log d} 
+ \frac{(\log\log d)-\ka}{40(\log d)^2} .
\end{align*}
\end{claim}
Now we are in a position to quickly deduce our bounds on $\alfm_d$.
\begin{proof}[Proof of Theorem~\ref{thm:alfm_bnd}]
According to \eqref{eq:Psi_g} we have 
\[ \psimain_d(\al_{d,a}) = \frac{a-\ka}{\log d} + g(x,y) . \]
Due to \eqref{eq:gxy_asymptotics}, for any $a \neq \ka$, the second term is negligible compared to the first term as $d \to \infty$. Thus $\psimain_d(\al_{d,a})$ is of order $1/ \log d$ and has the same sign as $a-\ka$. So (a) follows immediately from Lemma~\ref{lem:implications}(a). 

Now we turn to proving (b) and (c). Using the simple bound $\al \leq (2 \log d)/d$, we get the following for the expression that appears in the condition of Lemma~\ref{lem:implications}(b):
\begin{align}
\begin{split}
\frac 12 \al + \frac{d+1}{2}\al^2 
+ \frac 23 d \frac{\al^3}{1-2\al} 
&< \frac{\log d}{d} 
+ \frac{2 (\log d)^2}{d} 
+ \frac{2 (\log d)^2}{d^2}
+ \frac{8 (\log d)^3}{d^2} \\ 
&= \frac{(\log d)^2}{d} 
\left( \frac{1}{\log d} + 2 + \frac{2}{d} 
+ \frac{8 \log d}{d} \right) \\
&< \frac{2.2 (\log d)^2}{d} 
\quad \text{for } d \geq 1100 .
\label{eq:cond_b_bnd}
\end{split}
\end{align}
This is eventually smaller than $\psimain_d(\al_{d,a})$ if $a>\ka$, since the latter is of order $1/ \log d$. Therefore, we can apply Lemma~\ref{lem:implications}(b) and get $\alfm_d > \al_{d,a}$. 

As for $a=\ka$, by Claim~\ref{cl:g_pos} we have for $d \geq 2 \cdot 10^6$ that 
\[ \psimain_d(\al_{d,\ka}) 
> \frac{(\log\log d)-\ka}{40(\log d)^2} .\]
If this is larger than $2.2 (\log d)^2/d$, then we can apply Lemma~\ref{lem:implications}(b) to get $\alfm_d > \al_{d,\ka}$. So it remains to show that the following ratio is above $1$:
\[ \frac{(\log\log d)-\ka}{2.2 \cdot 40} 
\frac{d}{(\log d)^4} .\]
A direct check at $d=2 \cdot 10^6$ gives roughly $1.2$ (so it is comfortably above $1$), and the expression increases thereafter because both fractions are monotone increasing for $d > e^4$. \end{proof}
\begin{proof}[Proof of Claim~\ref{cl:g_pos}]
First note that for $a \geq \ka$ we have $y \leq x / \log d$, and hence 
\begin{align*}
x-y 
&\geq \left( 1-\frac{1}{\log d} \right) x ;\\
g(x,y) &> - y + \frac{(x-y)^2}{2} 
\geq -\frac{x}{\log d} 
+ \frac 12 \left( 1-\frac{1}{\log d} \right)^2 x^2 
= c_d \, \frac{x}{\log d} \text{, where } \\
c_d &\defeq -1 
+\frac 12 \left( 1-\frac{1}{\log d} \right)^2 x \log d \\
&= -1 + \frac 12 \left( 1-\frac{1}{\log d} \right)^2
\bigg( (\log\log d) - \ka \bigg) .
\end{align*}
Note that $c_d$ is monotone increasing and tends to $\infty$. It becomes positive between $10^6$ and $1.1 \cdot 10^6$. Therefore, 
\[ g(x,y)>0 
\text{, and hence }
\psimain_d( \al_{d,a} ) > \frac{a-\ka}{\log d}
\quad (d \geq 1.1 \cdot 10^6) .\]
Moreover, for $d \geq 2 \cdot 10^6$ we have $c_d > 0.026 > 1/40$. It follows that 
\[ \psimain_d(\al_{d,a}) > \frac{a-\ka}{\log d} 
+ \frac{(\log\log d)-\ka}{40(\log d)^2} . \]
\end{proof}

Finally, we deduce a specific lower bound for $\log( \al_{d,a} )$ that we used in Section~\ref{sec:technical}.
\begin{claim} \label{cl:log_al_minus}
Let $0 \leq a \leq \log\log d$. Then  
\begin{equation} \label{eq:log_al_minus} 
\log( \al_{d,a} ) \geq \log( \al_{d,\log\log d} ) 
> -\log d + \log\log d + \log(2) 
- \frac{\log\log d}{\log d}
- \frac{1/4}{\log d} .
\end{equation}
\end{claim}
\begin{proof}
First we show that for every $z>1$ 
\[ (\log z)^2 \leq \frac{4}{e^2} z .\]
Indeed, the derivative  
\[ \left( \frac{(\log z)^2}{z} \right)'= 
\frac{(\log z)(2- \log z)}{z^2} ,\]
vanishes only at $z=e^2$; so $(\log z)^2/z$ takes its maximum $4/e^2$ at $z=e^2$.  

In particular, for $z= \log d$ we get 
\begin{equation*} \label{eq:loglog_vs_log}
(\log\log d)^2 \leq \frac{4}{e^2} \log d ,
\end{equation*}
which, by \eqref{eq:log_low_bnd1}, implies that 
\begin{align*}
\log( 1-x+y ) &> \log(1-x) > -x - x^2 \\
&> \frac{- \log \log d + \ka}{\log d} 
- \frac{(\log\log d)^2}{(\log d)^2} \\
&\geq \frac{- \log \log d}{\log d} 
+ \frac{\ka - 4/e^2}{\log d} .
\end{align*} 
Since $\ka-4/e^2>-1/4$, \eqref{eq:log_al_minus} follows from \eqref{eq:log_al_a}. 
\end{proof}
\subsection{Connection to the Lambert W function}
It is worth mentioning that one can express the root of $\psimain_d$ using the \emph{Lambert $W$ function}, which, in turn, has an expansion in the form of a ``logarithmic power series''. 

To be more specific, let $W$ denote the inverse of the function $w \mapsto w \exp(w)$. Then it is easy to see that the root of $\psimain_d(\al)$ can be expressed as
\[ \frac{2}{d} W\left( \frac{e}{2} d \right) .\]
As we mentioned, $W(z)$ has an expansion: let 
$L_1 \defeq \log z$ and $L_2 \defeq \log\log z$; then 
\[ W(z) = L_1 - L_2 - \frac{L_2}{L_1} 
+ \frac{L_2(-2+L_2)}{2L_1^2} + \cdots .\]
Using this expansion, one may obtain arbitrarily precise estimates for the root of $\psimain_d$. 
%

\newpage

\section{Appendix: Factor of IID versus Typical} 
\label{sec:FIID_vs_typical}

As we mentioned in Section~\ref{sec:factors}, every FIID (factor of IID) process is typical. As for the reverse direction, there used to be a conjecture stating that every typical process on the $d$-regular tree $T_d$ is the weak limit of FIID processes. This can be phrased in various (essentially) equivalent ways. One such formulation uses local-global convergence, claiming that the sequence $(\Gb_{N,d})_N$ of random regular graphs a.s.~converges to the so-called Bernoulli graphing of $T_d$; see \cite[Section 7]{hatami2014limits}. Loosely speaking, these conjectures stated that whatever object is present in random $d$-regular graphs with high probability can be (approximately) obtained using a local algorithm. Specifically, Szegedy conjectured that, for each $d$, the density of factor of IID independent sets on $T_d$ may be arbitrarily close to $\alstar_d$; see the introduction of \cite{lyons2011perfect}. This was refuted by Gamarnik and Sudan for $d \geq d_0$ in \cite{gamarnik2017limits}. Note that the Gamarnik--Sudan result does not specify a threshold $d_0$, and whether the same is true for every $d$ (most notably for $d=3$) is still an open problem. Our explicit bounds on $\alstar_d$ enable us to settle the case $d \geq 403$.
\begin{theorem} \label{thm:fiid_ind_set}
If $d \geq 403$, then $T_d$ has a typical process that is not the weak limit of factor of IID processes: there exists $\al,\eps>0$ such that $T_d$ has a typical independent set of density $\al$ but every FIID independent set on $T_d$ has density at most $\al-\eps$.  
\end{theorem}

In this paragraph we say \emph{stable set} for an independent set in order to avoid the confusion with the probabilistic meaning of \emph{independent}. The key idea behind the Gamarnik--Sudan argument is the observation that a factor of IID stable set can be used to prove the existence of a pair of stable sets with a prescribed intersection density. Indeed, suppose that we have a factor that maps a collection of IID labels (placed on the vertices of $T_d$) to a random stable set of density $\al$. If we take two independent copies of the IID labels, then they are mapped to a pair of random stable sets that are independent, and hence the density of their intersection (i.e., the probability that a given vertex of $T_d$ is included in both stable sets) is $\al^2$. On the other hand, if we use the same underlying IID labels for both stable sets, then they, of course, coincide everywhere, and hence the intersection density is $\al$. In fact, it is easy to show that one can continuously interpolate between these two cases. More precisely, after we sample the first collection of IID labels, we define the second collection of IID labels as follows. We fix a probability $p \in [0,1]$. For each vertex, with probability $p$, we keep the original label, while with probability $1-p$ we re-sample the label (independently from everything else). Then the intersection density of the two stable sets is a continuous function of $p$ such that the function value is $\al^2$ at $p=0$ and $\al$ at $p=1$. It follows that for any $\be \in [\al^2,\al]$ we can choose $p$ such that the intersection density is $\be$. In particular, for every $\be \in [\al^2,\al]$ there is a typical process on $T_d$ that produces a pair of stable sets of density $\al$ with intersection density $\be$. It follows that $f_{d,\al}(\be)$ must be nonnegative; otherwise the expected number of such pairs in $\Gbconf_{N,d}$ would be exponentially small (see Section~\ref{sec:counting} for details).

\begin{figure}
\centering
\includegraphics[width=0.5\linewidth]{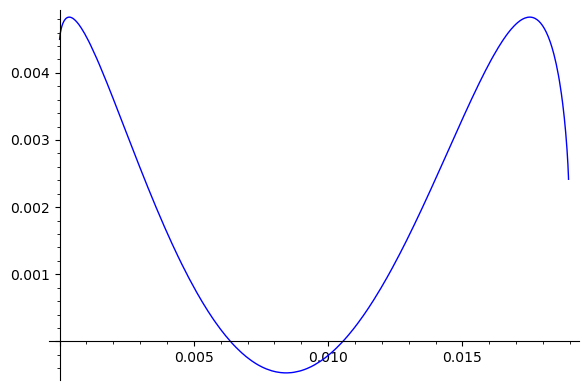}
\caption{The plot of $f_{d,\al}(\be)$ for $d=500$ and $\al=\alsm_d$. The two local maxima take the same value, and the local minimum is negative so the conditions of Corollary~\ref{cor:neg_min} are satsified.}
\label{fig:f_plot_d500}
\end{figure}

In conclusion, Gamarnik and Sudan essentially proved the following in \cite{gamarnik2017limits}.
\begin{lemma} \label{lem:fiid_ind_set}
Suppose that $T_d$ has a factor of IID independent set of density $\al$. Then 
\begin{equation*} 
f_{d,\al}(\be) \geq 0 
\quad \text{for all } \be \in [\al^2,\al] .
\end{equation*}
\end{lemma}
Combining this with Theorem~\ref{thm:sm_intro} of the current paper, we obtain the following.
\begin{corollary} \label{cor:neg_min}
Suppose that $d,\al$ are such that 
\begin{equation*} 
f_{d,\al} \text{ has two local maxima } 
\be_1^{\max}=\al^2 < \be_2^{\max}
\text{ and one local minimum } \be^{\min} .
\end{equation*}
If 
\begin{equation} \label{eq:both_cond}
f_{d,\al}\big( \be_1^{\max} \big) 
\geq f_{d,\al}\big( \be_2^{\max} \big)
\quad \text{and} \quad
f_{d,\al}\big( \be^{\min} \big) < 0 ,
\end{equation}
then $\alstar_d \geq \al$, while every FIID independent set has density at most $\al-\eps$ for some $\eps>0$.

It turns out that condition \eqref{eq:both_cond} holds for $\al=\alsm_d$ provided that $d \geq 458$. (See Figure~\ref{fig:f_plot_d500} for the plot of $f_{d,\al}$ in the case $d=500$, $\al=\alsm_d$.)
\end{corollary}

It follows that, above the threshold $d_0=458$, the maximal density of independent sets are larger for typical processes than for FIID processes. We get a somewhat better threshold if we use the augmented bound $\alaug_d$ (see Corollary~\ref{cor:alaug}). Setting $\al \defeq \alaug_d$, $f_{d,\al}\big( \be^{\min} \big)$ is negative starting at $d=403$ already, and hence Theorem~\ref{thm:fiid_ind_set} follows from Lemma~\ref{lem:fiid_ind_set}.

\subsection{Asymptotic FIID independence ratio}

As $d$ grows, the difference between typical and FIID independent sets gets substantial. Gamarnik and Sudan proved that the density of FIID independent sets is, asymptotically as $d\to \infty$, at most %
\[ \left( 1+\frac{1}{\sqrt{2}}+o_d(1) \right) 
\frac{\log d}{d} .\]
Essentially, they derived this from Lemma~\ref{lem:fiid_ind_set} via a lengthy computation. Below we include a streamlined version of this proof (see Claim~\ref{cl:gs_revisited} and its corollary). Later, Rahman and Vir\'ag \cite{rahman2017local} extended their method and proved the stronger (and sharp) bound 
\[ \big( 1+o_d(1) \big) \frac{\log d}{d} ,\]
which is roughly half of $\alstar_d$ (the maximal density of typical independent sets).

\begin{claim} \label{cl:gs_revisited}
Let $0<b<a$ be positive constants and set 
\[ \al = a \frac{\log d}{d} 
\quad \text{and} \quad 
\be = b \frac{\log d}{d} .\]
Then, as $d \to \infty$, we have 
\begin{equation*}
f_{d,\al}(\be) = \bigg( C(a,b) + o(1) \bigg) 
\frac{(\log d)^2}{d} 
\text{, where } 
C(a,b) \defeq \frac 12 - (a-1)^2 + \frac 12 (b-1)^2 .
\end{equation*}
\end{claim}
\begin{corollary}
Note that, for any fixed $a>1$, $C(a,b)$ takes its minimum at $b=1$. The minimum value 
\[ \min_{b \in [0,a]} C(a,b)=C(a,1)= \frac 12 - (a-1)^2 \]
is negative if and only if $a>1+\frac{1}{\sqrt 2}$. For any such $a$, Lemma~\ref{lem:fiid_ind_set} can be applied (for sufficiently large $d$), showing that a FIID independent set on $T_d$ has density at most 
\[ \bigg( 1+ \frac{1}{\sqrt{2}} +o(1) \bigg) 
\frac{\log d}{d} .\]
\end{corollary}
\begin{proof}[Proof of Claim~\ref{cl:gs_revisited}]
In this proof $\approx$ means that the equality holds with an error $o\big( (\log d)^2 / d \big)$ as $d \to \infty$. With this notation we need to show that $f_{d,\al}(\be) \approx C(a,b) (\log d)^2/d$. Note that, given any constant $x$, 
\[ \text{for }
\xi = x \frac{\log d}{d} : \quad
h(\xi) = x \frac{\log d}{d} 
\big( \log d - \log\log d - \log x \big) 
\approx x \frac{(\log d)^2}{d} .\]
Recall that $f_{d,\al}(\be)= \fhat_{d,\al}(\be, \Ga(\be))$; see \eqref{eq:f_al}. So let $\ga \defeq \Ga(\be)$. From \eqref{eq:Ga_second} we get 
\[ \ga 
= \big( 1+o(1) \big) (\al-\be)^2 
= \big( 1+o(1) \big) (a-b)^2 \frac{(\log d)^2}{d^2} .\]
From \eqref{eq:fhat} we have 
\begin{align*}
f_{d,\al}(\be) 
&= h(\be)+2h(\al-\be)+h(1-2\al+\be) + d h(\ga) \\
&+ 2d\bigg( h(\al-\be-\ga) - h(\al-\be) \bigg) \\
&+ \frac{d}{2} h(1-4\al+2\be+2\ga) - dh(1-2\al+\be) .
\end{align*}
Now we estimate the various terms of the right-hand side: 
\begin{align*}
& h(\be)+2h(\al-\be) \approx (2a-b) \frac{(\log d)^2}{d} ;\\
& h(1-2\al+\be) \approx 2\al-\be \approx 0 ;\\
& d h(\ga) \approx -2d \ga \log(\al-\be) ;\\
& 2d\bigg( h(\al-\be-\ga) - h(\al-\be) \bigg)
\approx 2d\ga \big( 1 + \log(\al-\be) \big) ;\\
& \frac{d}{2} h(1-4\al+2\be+2\ga) - dh(1-2\al+\be) \\ 
&\quad \approx 
\frac{d}{2} (4\al-2\be-2\ga) 
-\frac{d}{4} (4\al-2\be-2\ga)^2 
-d (2\al-\be)
+ \frac{d}{2} (2\al-\be)^2 \\
&\quad \approx 
-d\ga - \frac{d}{2} (2\al-\be)^2 
= -d\ga - \frac{1}{2} (2a-b)^2 \frac{(\log d)^2}{d} .
\end{align*} 
Putting everything together we obtain 
\begin{align*}
f_{d,\al}(\be) &\approx 
(2a-b) \frac{(\log d)^2}{d} + d\ga 
- \frac{1}{2} (2a-b)^2 \frac{(\log d)^2}{d} \\
&\approx \bigg( 2a-b + (a-b)^2 - \frac{(2a-b)^2}{2} \bigg) 
\frac{(\log d)^2}{d} \text{, where}
\end{align*}
\[ 2a-b + (a-b)^2 - \frac{(2a-b)^2}{2} 
= -a^2+2a +\frac 12 b^2-b 
= \frac 12 - (a-1)^2 + \frac 12 (b-1)^2 ,\]
and the proof is complete. 
\end{proof}

\end{document}